\input amstex
\documentstyle{amsppt}
\magnification=1150
\voffset=-.2in
%\vsize=7.5in
%\pagewidth{6.5in}
%\pageheight{8in}
\input xy
\xyoption{v2}
\NoBlackBoxes
\NoRunningHeads

%%%%%%%%%%%%%%%%%%%%%%%%%%%%%%%%%%%%%%%%%%%%%%%%%%%%%%%

%macros for text, and abbreviations of other macros used below.

                               %

\define\wt{\widetilde}                                %
\define\ov{\overline}                                 %
\define\tagit#1{\tag "(#1)"}
\define\spec#1{\underline{\underline{#1}}}

\define\Wit2#1{W_*(#1) \tsize{[\frac 12]}}
\define\Wt2#1{\underline{\underline{W}}(#1)\tsize{[\frac12]}}

\define\prf{\demo{\underbar{Proof}}}
\define\endpf{\enddemo}
\define\dfn#1{\definition{\bf\underbar{Definition #1}}}
\define\enddef{\enddefinition}
\define\BbbX{\Bbb X}                                  %
\define\gm{\Gamma}                                    %
\define\surj{\twoheadrightarrow}                      %
\define\inj{\rightarrowtail}                          %

\document
\baselineskip 20pt
\hsize = 6.5 true in

%%%%%%%%%%%%%%%%%%%%%%%%%%%%%%%%%%%%%%%%%%%%%%%%%

\topmatter
\title
Polynomially bounded cohomology and discrete groups
\endtitle
\vskip.2in
\author
C. Ogle
\endauthor
\affil
The Ohio State University\\
\\
May 2001; revised March, 2004
\endaffil
\abstract
We establish the homological foundations for studying polynomially bounded group cohomology, and show that the natural map from $PH^*(G;\Bbb Q)$ to $H^*(G;\Bbb Q)$ is an isomorphism for a certain class of groups.
\endabstract
\toc
\widestnumber\subhead{3.2.1}
\specialhead {} Introduction
\endspecialhead
\head 1. Polynomially bounded group cohomology
\endhead
\subhead 1.1. Basic results in polynomially bounded cohomology and the 
Leray-Serre spectral sequence
\endsubhead
\subhead 1.2. The obstruction to injectivity
\endsubhead
\subhead 1.3. Analysis of the obstruction
\endsubhead
\subhead 1.4. The map $\alpha_1(G,L_G)$
\endsubhead
\head 2. Higher Dehn functions
\endhead
\subhead 2.1. Dehn functions and simplicial resolutions
\endsubhead
\subhead 2.2. Higher Dehn functions and cohomology
\endsubhead
\subhead 2.3. Linearly and uniformly bounded cohomology
\endsubhead
\head Appendix. \hskip.5in Type $P$ resolutions
\endhead
\endtoc
\email
ogle@math.mps.ohio-state.edu
\endemail
\endtopmatter

%%%%%%%%%%%%%%%%%%%%%%%%%%%%%%%%%%%%%%%%%%%%%%%%%
%%%%%%%%%%%%%%%%%%%%%%%%%%%%%%%%%%%%%%%%%%%%%%%%%

\newpage
\vskip.2in

\centerline{\bf\underbar{Introduction}}
\vskip.3in

The cohomology of a discrete group $G$ with coefficients in a $G$-module $A$ can be defined
in various equivalent ways. Typically one first constructs a cocomplex, which for now we
will label $(C^*(G;A),\delta^*)$; the cohomology of $G$ with coefficients in $A$ is then the
cohomology of this complex.

Suppose $G$ is a countable group equipped with word-length function $L$.
Given the pair $(G,L)$, one can consider various refinements of
this cocomplex which involve a growth condition on the level of cochains. The most restrictive
is a uniform bound. This condition defines a subcomplex of bounded cochains which is already
quite interesting and has been extensively studied over the last thirty years ([Gr], [Gr1], [I1], [N1], [P1], [Mi1], [Mi2]). Less restrictive (and also less studied) is the
case when the growth rate on the level of cochains is \underbar{polynomial}. This growth
condition is related to the Novikov conjecture, as shown in [CM].
For suitable $A$ (as defined in section 1.1), one has a
natural subcocomplex $PC_{L}^*(G;A)\subset C^*(G;A)$ consisting of cochains of polynomial growth
with respect to $L$,
and the inclusion is functorial with respect to polynomially bounded group homomorphisms in the
first coordinate, and polynomially bounded module homomorphisms in the second coordinate.
The resulting cohomology groups of $PC_{L}^*(G;A)$ are denoted by $PH_{L}^*(G;A)$;
in general the cocomplex $PC_{L}^*(G;A)$ and therefore also its cohomology groups
depend on the choice of word-length function $L$. The inclusion of
cocomplexes induces a transformation $\eta(G,L;A)^* : PH_{L}^*(G;A)\to H^*(G;A)$.
Of most interest to us is the case when $A = \Bbb Q$. We consider three sucessively weaker
conditions one could ask of the group $G$:
\vskip.1in

\noindent\underbar{(PC1)} The map $\eta(G,L;\Bbb Q)^*$ is an isomorphism.
\vskip.05in

\noindent\underbar{(PC2)} The map $\eta(G,L;\Bbb Q)^*$ is an epimorphism.
\vskip.05in

\noindent\underbar{(PC3)} For every $0\ne x\in H_*(G;\Bbb Q)$ there is a 
$y\in PH_L^*(G;\Bbb Q)$ with
$<\eta(G,L;\Bbb Q)(y), x>\ne 0$.
\vskip.1in

The
dual of $\eta(G,L;\Bbb Q)^*$ is a $\Bbb Q$-vector space map
$(\eta(G,L;\Bbb Q)^*)^* : (H^*(G;\Bbb Q))^*\to (PH_L^*(G;\Bbb Q))^*$, and we can form the
composition
$$
\alpha_*(G,L) : H_*(G;\Bbb Q)\inj (H^*(G;\Bbb Q))^*\to (PH_L^*(G;\Bbb Q))^*
\tag0.1
$$
where the first map is induced by the Universal Coefficient Theorem, and the second is
$(\eta(G,L;\Bbb Q)^*)^*$. The above 3 conditions can be rephrased as
\vskip.1in

\underbar{(PC1)} $(\eta(G,L;\Bbb Q)^*)^*$ is an isomorphism.
\vskip.05in

\underbar{(PC2)} $(\eta(G,L;\Bbb Q)^*)^*$ is a monomorphism.
\vskip.05in

\underbar{(PC3)} $\alpha_*(G,L)$ is a monomorphism.
\vskip.1in

For certain geometric groups it is feasible
to verify property (PC1), which we do in this paper. A weaker
condition is (PC2); this is the condition (PC) of [CM] and equivalent to (PC3) when the rational
homology groups of $G$ are degreewise finitely generated. However when the rational homology of $G$ is not  finitely generated in each degree, (PC2) is more restrictive than (PC3). For example, if $G$ is a free group on a countably infinite set of generators and $L$ a word-length metric on $G$ (see below), then (PC3) holds but (PC2) fails. Also, injectivity of the map in (PC3) is sensitive to the choice of word-length, and injectivity
may hold for some choices of word-length but not for others. To illustrate, we see that the
condition is obviously satisfied for $G = \Bbb Z$ with the standard word-length. However, if we
use instead a word-length which depends logarithmically on the standard one, then with respect
to this word-length $PH^1(\Bbb Z;\Bbb Q) = 0$ and so (PC3) fails.
The issue of injectivity of $\alpha_i(G,L)$ is related to the
\underbar{Dehn function} $f_G$ of $G$. This function, introduced and studied by Gersten
([Ger1], [Ger2]) is defined in terms of the presentation of the group. Given a word whose image
in $G$ is trivial, the Dehn function measures the increase in word-length when one writes this
word as a minimal product of conjugates of relators occuring in the relator set of the
presentation. Although the Dehn function itself depends on the presenation, the linear
equivalence class in which it lies does not [Ger1]. Thus up to such equivalence, one may simply
refer to the Dehn function of $G$. The word-problem for $G$ is solvable iff $G$ has a
recursively enumerable Dehn function, and solvable in polynomial time iff $f_G$ is polynomial.
All known computations support the following conjecture

\proclaim{\bf\underbar{Conjecture A}} If $\Cal P = <\Cal S | \Cal W>$ is a finite presentation
of $G$ with polynomial Dehn function $f_G$, then $\alpha_i(G,L^{st}_G)$ is an injection for all $i\ge 0$.
\endproclaim

In this paper we establish a framework for proving Conjecture A.
First, in section 1.1 we establish some basic
results in p-bounded homological algebra. Primarily, we construct the Serre
spectral sequence associated to a short-exact sequence of groups with word-length (as defined
in that section); the existence of the proper $E^1_{pq}$-term for $q >1$ requires an additional hypothesis, but for $q=0,1$ the spectral sequence takes the usual form, which leads to a five-term exact sequence analogous to the one in ordinary cohomology (cf. [N1], [P1] for the corresponding spectral sequence in bounded cohomology). We also prove a Comparison Theorem, which tells us under what conditions a resolution can be used to compute p-bounded cohomology. 
Section 1.2 uses the five-term sequence to identify the obstruction to injectivity of $\alpha_i(G,L)$ in even dimensions. In section 1.3, we show that for groups with polynomial Dehn function, a related obstruction vanishes. The results of this section are in preliminary form; a detailed account will appear in a sequel to this paper. In section 1.4 we verify the injectivity of $\alpha_1(G,L)$ for a suitable choice of $L$ when the $H_1(G)$ is finitely-generated.

Section 2 contains various results related to Dehn functions. In section 2.1, we show how type $P$ resolutions (Appendix A) can be used to define the higher Dehn functions for groups of type $FP^{\infty}$. The constructions in this section are then used in section 2.2 to prove 

\proclaim{\bf\underbar{Theorem B}} If $\Cal P = <\Cal S | \Cal W>$ is a finite presentation
of $G$ with polynomial Dehn
function $f_G$, then $\eta(G,L;A)^i$ is an isomorphism for any p-semi-normed (p.s\.) $G$-module $A$ (defined in section 1.1) and $i = 1,2$. 
\endproclaim

In general $\eta(G,L;\Bbb Q)^2$ fails to be surjective when $f_G$ is non-polynomial [Ger3]. In fact, $\alpha_2(G,L^{st}_G)$ is not injective for the example Gersten constructs in that paper.
In section 2.2, using [G5] we show

\proclaim{\bf\underbar{Theorem C}} If $G$ admits a bounded combing (in particular, if
$G$ is automatic), then $\eta(G,L;A)^*$ is an isomorphism for all p.s\. $G$-modules $A$.
\endproclaim

This map is also an isomorphism when $G$ is nilpotent.
In section 2.3, we define linearly bounded (or Lipschitz) cohomology $LH^*(G;A)$ for appropriate coefficient modules $A$. As with p-bounded cohomology the inclusion map on the cochain level induces a natural homomorphism
$$
\eta_{lin}(G,L;A)^*: LH^*(G;A)\to H^*(G;A)
$$
It is a theorem due to Gromov that $f_G$ is linear iff $G$ is
word-hyperbolic. Recently a complete cohomological characterization of
word-hyperbolic groups has  been obtained by I. Mineyev in [Mi1]. Using the result of [Mi3] we show

\proclaim{\bf\underbar{Theorem D}} If $G$ is word-hyperbolic, then $\eta_{lin}(G,L;A)^*$ is an isomorphism for all l.s\. $G$-modules $A$.
\endproclaim
In fact, the result in [Mi1] suggest the stronger statement that $\eta_{lin}(G,L;A)^*$ is an isomorphism for all $A$ iff $G$ is word-hyperbolic.

In the appendix we cover the definition and formal properties of type $P$ resolutions as developed in [O1].

A remark on notation: throughout the paper we write $PC^*_{L}(\;\;\;;A)$ resp.
$PH^*_{L}(\;\;\;;A)$ as $PC^*(\;\;\;;A)$ resp. $PH^*(\;\;\;;A)$, and $\alpha(G,L)$ as
$\alpha(G)$ unless we need to emphasize a particular word-length function.

%%%%%%%%%%%%%%%%%%%%%%%%%%%%%%%%%%%%%%%%%%%%%%%%%
%%%%%%%%%%%%%%%%%%%%%%%%%%%%%%%%%%%%%%%%%%%%%%%%%

\newpage
\vskip.2in

\centerline{\bf\underbar{Polynomially bounded group cohomology}}
\vskip.3in

\subhead
1.1 Basic results in polynomially bounded cohomology and the Leray-Serre spectral sequence
\endsubhead
\vskip.2in

If $S$ is a generating set for a free group $F$ and
$f : S\to \Bbb N^+$ a function, then $S$ and $f$ determine a word-length function $L_F$
on $F$ given by
$$
\gather
L_F(id) = 0\\
L_F(x) = f(x)\quad \text{if}\ x\ \text{or}\ x^{-1}\ \text{is in}\ S\\
L_F(g) = \sum_{i=1}^r f(x_i)\qquad\text{where}\ x_1x_2\dots x_r\ \text{is the unique reduced word representing
$g$}
\endgather
$$
such a word-length function on $F$ is referred to as a \underbar{word-length metric}. 
If $F'\subset F$ is a subgroup of $F$ equipped with a word-length metric $L_F$, then the restriction of $L_F$ to $F'$ defines an \underbar{induced metric} $L_{F'}$ on $F'$. Finally, if $p:F'\surj G$ is a surjection of $F'$ to $G$, 
then $L_{F'}$ determines a word-length function $L_G$ on $G$ by $L_G(g) = min\{L_{F'}(f)\ |\ p(f) = g\}$.
Any non-degenerate word-length function on $G$ may be realized in this fashion for an appropriate choice of $F$, $f$ and $p$. When the set $S$ is finite and $f(x) = 1$ for each $x\in S$, $L_G$ is referred to as the
standard word-length function $L^{st}_G$ associated with the set of generators $S$.
Note that $L_G$ depends only on the pair $(S,f)$, so that if $<S|W>$ and $<S|W'>$ are two
presentations of $G$ which have the same set of generators and weight function $f$, then
the induced word-length functions will also be the same.

We will use the notation $A[S]$ to denote the free $A$-module with
basis $S$ for a countable set $S$.
In particular, if $\Bbb Q[G]$ is the rational group algebra of $G$, $\Bbb Q[G][S]$ is the free
$\Bbb Q[G]$-module with basis $S$.

A \underbar{weighted set} is a pair $(S,f_S)$ where $f_S : S\to \Bbb R^+$
is a function, referred to as the weight function. When the weight
function is understood, we write $(S,f_S)$ simply as $S$.

A homomorphism $(G,L)\to (G',L')$ of groups with word-length will mean a homomorphism
$f : G\to G'$ for which $L'(f(g)) = min\{L(h) \ |\ f(h) = f(g)\}$. Thus, if $f$ is a
monomorphism it preserves word-length, and if $f$ is an epimorphism, $L'$ is the
word-length function induced by $f$ and $L$. A short-exact sequence of groups with
word-length is a sequence of morhisms of groups with word-length
$$
(K,L_K)\inj (G,L_G)\surj (N,L_N)
$$
where the underlying sequence of groups and group homomorphisms is
short-exact. Note that if $K\inj G\surj N$ is a short-exact sequence of groups
and $L_G$ a word-length function on $G$, then there exist unique word-length functions
$L_K$ resp. $L_N$ on $K$ resp. $N$ making
$(K,L_K)\inj (G,L_G)\surj (N,L_N)$ a short-exact sequence of groups with word-length.

A \underbar{semi-norm} $\eta$ on a $k$ vector space $V$
($k\subset \Bbb R$) is a map $\eta : V\to \Bbb R_+$ satisfying i)
$\eta(a+b)\le \eta(a) + \eta(b)$ and ii) $\eta(\lambda a)\le
|\lambda|\eta(a)$ for all $a,b\in V$ and $\lambda\in \Bbb Q$.
\vskip.1in

Before proceeding, we illustrate an essential homological difference between the notions of \lq\lq bounded\rq\rq and \lq\lq p-bounded\rq\rq. If $(S,f_S)$ is a weighted set and $(V,\|\,\,\,\|)$ a normed vector space, one may define $BHom(S,V)$ the set of bounded morphisms from $S$ to $V$, and a larger space $PHom(S,V)$, the set of p-bounded morphisms from $S$ to $V$. $\phi:S\to V$ is p-bounded if there is a polynomial $p$ such that $\|\phi(s)\|\le p(f_S(s))$ for all $s\in S$. Then $\phi$ is bounded if we can take the polynomial to be a constant function. So $BHom(S,V)$ is again a normed vector space with respect  to the sup norm, allowing one to construct spaces such as $BHom(S', BHom(S,V))$. The larger space $PHom(S,V)$ has no natural norm. However, it does have an obvious collection of semi-norms given by $\eta_s(\phi) = \|\phi(s)\|$. This suggests that in the p-bounded setting, one needs to work with semi-normed modules of a sufficiently general type in order to define iterated $Hom$ spaces such as $PHom(S', PHom(S,V))$ ( a necessary construction for the developement of the Serre spectral sequence in p-bounded cohomology). 

\vskip.2in

\dfn{1.1.1} A \underbar{p-semi-normed $G$-module}, or p.s\.
$G$-module, is a $\Bbb Q[G]$-module $M$ equipped with a collection
of semi-norms $\{\eta_x\}_{x\in \Lambda}$ indexed on a countable
$G$-set $\Lambda = \Lambda_M$. The semi-norms satisfy the
following properties \vskip.01in

i) If $\eta_{x_1},\dots,\eta_{x_n}\in \Lambda$ and $\lambda_1,\dots,\lambda_n\in
\Bbb N^+$, then there is a constant $C$ and an $\eta_y\in\Lambda$ with
$$
C\eta_y\ge \lambda_1\eta_{x_1}+\dots +\lambda_n\eta_{x_n}
$$
\vskip.05in

ii) there exist constants $C,n > 0$ such that for all
$x'\in\Lambda$ there is an $x\in \Lambda$
with
$$
\eta_{gx'}(hm)\le C(1+L(h))^n\eta_{ghx}(m)
$$
for all $g,h\in G$ and $m\in M$, where $L$ is the word-length function on $G$.
\enddef
\vskip.05in

When $G = \{id\}$ we will refer to $M$ as a p.s\. module. Any p.s\. $G$-module
is a p.s\. module by forgetting the $G$-module structure.

\dfn{1.1.2}
A \underbar{homomorphism $f : M\to M'$ of p.s\. $G$-modules} is a $\Bbb Q[G]$-module
homomorphism which is
\underbar{p-bounded}; i.e. there exists  $C_1,C_2$ and $n > 0$ such that for all
$x'\in \Lambda_{M'}$ there exists $x\in \Lambda_M$ with
$$
\eta_{gx'}(f(a))\le C_1C_2(1+\eta_{gx}(a))^n(1+L(g))^n
$$
for all $g\in G$ and $a\in M$. In this inequality, the constants $C_1$ and $n$ may vary with $f$
but are independent of $x'$, while $C_2$ depends only on $x'$.
\enddef
\vskip.05in

Two p.s\. $G$-modules $M$, $M'$are \underbar{isomorphic}
if there exist homomorphisms of p.s\. $G$-modules $f:M\to M'$, $f' : M'\to M$ with
$f\circ f' = id_{M'}$, $f'\circ f = id_M$. By (1.1.1) i), the collection of all p.s\. $G$-module
homomorphisms from $M$ to $M'$ forms a vector space over $\Bbb Q$ which
we denote $PHom_G(M,M')$. Dropping the requirement that the maps be $G$-equivariant, we
get the $\Bbb Q$-vector space of p-bounded maps from $M$ to $M'$ which we denote
simply as $PHom(M,M')$. The same conventions apply for $Hom$ in place of $PHom$.
If $M'$ is a sub-$G$-module of $M$ and $M'' = M/M'$, we may define a p.s\. $G$-module
structure on $M''$ by setting $\Lambda_{M''} =\Lambda_M$ and for all $x\in \Lambda_{M''}$
defining
$$
\eta_x(\ov m) = min\{\eta_x(m)\,|\, p(m) = \ov m\}
\tag1.1.3
$$
where $p : M\surj M/M'$ is the projection. The reader may verify that this defines a
p.s\. $G$-module structure on $M''$. We refer to this as the
\underbar{quotient p.s\. $G$-module
structure} induced by $M$ and the projection $p$.

It will sometimes be the case that a $G$-module $M$ comes equipped with two p.s\. $G$-module structures, which we may denote $\Cal S_1$ and $\Cal S_2$. Let $M_{\Cal S_i}$ denote $M$ equipped with structure $\Cal S_i$. We say that the two structures are \underbar{equivalent} if the identity map on $M$ induces an isomorphism of p.s\. $G$-modules $M_{\Cal S_1}\overset{id}\to\to M_{\Cal S_2}$.

\dfn{1.1.4}
A \underbar{free p-bounded $G$-module}, or p.f\. $G$-module, is a free $\Bbb Q[G]$-module
$P = \Bbb Q[G][S]$ with countable basis $S\ne \emptyset$ equipped with a weight function
$w_S : S\to \Bbb R^+$. The indexing set for the semi-norms on $P$ is
$\Lambda_P = \{*\}$ equipped with trivial $G$ action.
The unique semi-norm on $P$ is
$$
\left |\sum \lambda_i g_i s_i\right | = \sum |\lambda_i|(1 + |g_is_i|)
$$
where $|gs| = L(g) + w_S(s)$.
\enddef
\vskip.05in

 In the
special case $S = \{*\}$,  we adopt the convention that $|*| = 1$. In particular, when $S = \{*\}$ and $G = \{id\}$, this defines a semi-norm on $P = \Bbb Q$ given by
$|q|_* = 2|q|$. A p.f\. $G$-module is a p.s\. $G$-module. To see this, note that (1.1.1) i) is trivially satisfied because the indexing set has only one element, and (1.1.1) ii) follows from the inequality
$$
\gather
\left |h\left(\sum \lambda_i g_i s_i\right)\right | = \sum |\lambda_i|(1 + |hg_is_i|)\\
\le (1 + L(h))\sum |\lambda_i|(1 + |g_is_i|)
\endgather
$$
Suppose now that $P = \Bbb Q[G][S]$ is a p.f\. $G$-module and
$M$ a p.s\. $G$-module.
Let $\Lambda_M$ be the indexing set for the semi-norms on $M$ and let
$T = \{gs\ |\ g\in G,s\in S\}$. Associated to a finite subset $U\subseteq T$ and an element
$x\in\Lambda_M$ is the semi-norm on $PHom(P,M)$ given by
$$
\eta_{(x,U)}(f) = \sum_{t\in U}\eta_x(f(t))
$$
The indexing set for this collection of semi-norms is $\Lambda_M\times \Cal P(T)$,
where $\Cal P(T)$ denotes the set of finite subsets of $T$. The $G$-action on $PHom(P,M)$
is given by $g\cdot f(p) = gf(g^{-1}p)$, and the $G$ action on the index set
$\Lambda_M\times \Cal P(T)$
is given by $g(x,U) = (gx,g^{-1}U)$.

\proclaim{\bf\underbar{Proposition 1.1.5}} The above defines a p.s\. $G$-module structure
on $PHom(P,M)$.
\endproclaim

\prf For
$g,h\in G$, the series of inequalities
$$
\gather
\eta_{h(x',U)}(g\cdot f)\\
= \eta_{(hx',h^{-1}U)}(g\cdot f)\\
= \sum_{t_i\in U} \eta_{hx'}(gf(g^{-1}h^{-1}t_i))\\
\le C(1 + L(g))^{n}\left(\sum_{t_i\in U} \eta_{hgx}(f(g^{-1}h^{-1}t_i))\right)\\
= C(1 + L(g))^{n}\eta_{(hgx,g^{-1}h^{-1}U)}(f)\\
= C(1 + L(g))^{n}\eta_{hg(x,U)}(f)
\endgather
$$
implies the $G$-action is p-bounded in the sense of Def. 1.1.1. ii).
To verify (1.1.1) i) we suppose given numbers $\lambda_i > 0$ and semi-norms
$\eta_{(x_j,U_j)}$, $1\le j\le N$. Let $U = \cup_j U_j$ and
choose $y\in\Lambda_M$ and a constant $C'$ with $\sum_{j=1}^N  \lambda_j\eta_{x_j}\le C\eta_y$.
Letting $C = NC'$, one has the inequality
$$
\sum_{j=1}^N \lambda_j\eta_{(x_j,U_j)} \le C\eta_{(y,U)}
$$
\vskip-.35in\hfill //
\endpf
\vskip.35in

\proclaim{\bf\underbar{Proposition 1.1.6}} If $P = \Bbb Q[G][S]$ is a  p.f\. $G$-module and
$M$ a p.s\. $G$-module, then there is an isomorphism of vector spaces over $\Bbb Q$
$$
PHom_G(P,M)\cong PHom(\Bbb Q[S],M)
$$
\endproclaim

\prf This is the p-bounded analogue of a standard fact from homological algebra. The map from
left to right is given by restriction to the subspace $\Bbb Q[S]$; this restriction map is
obviously p-bounded. The map in the other direction is the inflation map. A
p-bounded map $f : \Bbb Q[S]\to M$ defines a $G$-module homomorphism $\wt f : P\to M$
given on basis elements by $\wt f (gs) = gf(s)$. The inequalities
$$
\gather
\eta_{hx'}\left(\wt f\left(\sum\lambda_i g_i s_i\right)\right) =
\eta_{hx'}\left(\sum\lambda_i g_i f(s_i)\right)\\
\le \sum |\lambda_i|C(1 + L(g_i))^n\eta_{hg_ix}(f(s_i))\\
\le \sum |\lambda_i|C(1 + L(g_i))^nC_1C_2(1 + |s_i|)^{n'}(1+L(h))^{n'}(1+L(g_i))^{n'}\\
\le (CC_1)C_2\left(1+\left |\sum \lambda_i g_i s_i\right |\right)^{n+n'}(1+L(h))^{n+n'}
\endgather
$$
imply $\wt f$ is p-bounded. \hfill //
\endpf

{\bf\underbar{Notes}}: i) By the last proposition, the p.s\. module structure on
$PHom(\Bbb Q[S],M)$  determines one on $PHom_G(P,M)$. There is also the induced p.s\. module
structure coming from the inclusion of $PHom_G(P,M)$ into $PHom(P,M)$ as the fixed-point set
under the action of $G$. The indexing set for the first is $\Lambda_M\times\Cal P(S)$, which
includes into that of the second, which is $\Lambda_M\times\Cal P(T)$. It is an
easy exercise to verify that these two structures are equivalent in the sense defined above.
\vskip.05in
ii) Inspection of the proofs of the previous two propositions show that $f:P\to M$ is p-bounded
precisely when it is p-bounded on the weighted set $T$, where $T$ is as in Proposition 1.1.5.
If $f$ is $G$-equivariant, then it lies in $PHom_G(P,M)$ precisely when it is a $G$-module
homomorphism which is p-bounded on the weighted set $S$.

We next discuss short-exact sequences.

\dfn{1.1.7} An \underbar{admissible monomorphism} $i : M'\inj M$
of p.s\. $G$-modules is a $G$-module monomorphism where
\vskip.03in i) $\Lambda_{M'} = \Lambda_M$; \vskip.03in ii) the
semi-norm $\eta_x$ on $M'$ is given by the restriction of $\eta_x$
on $M$ to $im(i)$. \vskip.05in An \underbar{admissible
epimorphism} $M\surj M''$ is an epimorphism which \vskip.03in i)
is a p.s. $G$-module homomorphism, and \vskip.03in ii)  admits a
section of p.s\. modules (i.e., a p-bounded homomorphism which is
not necessarily equivariant).
\enddef

In particular, the semi-norms on $M''$ may be given separately and are not
necessarily induced by the semi-norms on $M$.
A \underbar{short-exact sequence} of p.s\. $G$-modules
$$
M'\overset i\to\inj M\overset j\to\surj M''
$$
is then a short-exact sequence of $\Bbb Q[G]$-modules consisting of an admissible monomorphism
followed by an admissible epimorphism.

\proclaim{\bf\underbar{Lemma 1.1.8}} If $P$ is a p.f\. $G$-module and
$M'\overset i\to\inj M\overset j\to\surj M''$ a
short-exact sequence of p.s\. $G$-modules, then $i$ and $j$ induce a short-exact sequence
$$
PHom_G(P,M')\overset i_*\to\inj PHom_G(P,M)\overset j_*\to\surj PHom_G(P,M'')
$$
of $\Bbb Q$-modules.
\endproclaim

\prf Write $P$ as $\Bbb Q[G][S]$. By Proposition 1.1.6, the above sequence is isomorphic to
the sequence
$$
PHom(\Bbb Q[S],M')\overset i_*\to\inj PHom(\Bbb Q[S],M)\overset j_*\to\surj PHom(\Bbb Q[S],M'')
$$
obviously the first map is injective. The existence of a p-bounded section from $M''$ to $M$
implies the surjectivity of the second map. Lastly, if $f\in PHom(\Bbb Q[S],M)$ maps to zero
in $PHom(\Bbb Q[S],M'')$, its image lies in $i(M')$. Because the semi-norms on $M'$ are induced
by those on $M$ via the inclusion, the unique map $f' : \Bbb Q[S]\to M'$ for which
$f = i\circ f'$ is also p-bounded, and therefore an element of $PHom(\Bbb
Q[S],M')$. \hfill //
\endpf

[\underbar{Addendum to Lemma 1.1.8}: Although it will not be needed for what
follows, we note that the short-exact sequence in the above Lemma is
actually a short-exact sequence of p.s\. modules, where the p.s\. module
structure of each term is that described in the note following Proposition 1.1.6.]
\vskip.05in

A \underbar{p.s\. $G$-complex} is a $\Bbb Q[G]$-complex
$M_* = (M_*,d_*)$ where each $M_n$ is
a p.s. $G$-module and each boundary map $d_n: M_n\to M_{n-1}$ is a p.s\.
$G$-module homomorphism. A
p.s\. $G$-cocomplex $M^* = (M^*,\delta^*)$ is defined in exactly the same manner, with
$\delta^n : M^n\to M^{n+1}$. Given a p.s\.
$G$-complex $M_* = (M_*,d_*)$ and
a p.s\. $G$-module $M'$, we have a well-defined cocomplex
$$
PHom_G(M_*,M')
$$
with corresponding cohomology groups $PH^*_G(M_*;M')$, which are the p-bounded
$G$-equivariant cohomology groups of $M_*$ with coefficients in $M'$.

\dfn{1.1.9} A \underbar{p.f\. resolution} of $\Bbb Q$ over $\Bbb Q[G]$ is
a resolution $(R_*,d_*)$ of $\Bbb Q$ over $\Bbb Q[G]$ where each $R_n$ is a p.f\. $G$-module
and each $d_n : R_n\to R_{n-1}$ a p.s\. $G$-module homomorphism. In addition, we require that
$(R_*,d_*)$ admits a p-bounded chain contraction
$s_* = \{s_n : R_n\to R_{n+1}\}_{n\ge 0}$ as a p.s\. complex.
\enddef

The standard
non-homogeneous bar resolution over $\Bbb Q$, which we write as $(EG)_* = C_*(EG.;\Bbb Q)$,
provides an example of such a resolution. Precisely, for each $n\ge
0$ we identify $C_n(EG.;\Bbb Q)$ as the free p.s\. $G$-module on the set
$S_n$
$$
C_n(EG.;\Bbb Q) = \Bbb Q[G][S_n]
$$
where $S_n = \{(1,g_1,g_2,\dots g_n)\in EG_n = (G)^{n+1}\}$
and $G$ acts by left multiplication in the left-most coordinate:
$$
g(g_0,g_1,g_2,\dots,g_n) = (gg_0,g_1,g_2,\dots,g_n)
$$
the weight function on $S_n$ is given by
$$
f_{S_n}((1,g_1,g_2,\dots,g_n)) = 1 + \sum_{i=1}^n L_G(g_i)
$$
The differential $d_n$ defined on basis elements by
$$
d_n(g_0,g_1,g_2,\dots g_n) =
\left(\sum_{i=0}^{n-1}(-1)^i(g_0,\dots,g_ig_{i+1},\dots,g_n)\right) +
(-1)^n(g_0,g_1,\dots,g_{n-1})
$$
is p-bounded for each $n$, and $G$-equivariant. The standard section $s_*$ is
defined on basis elements by
$$
s_n((g_0,g_1,\dots ,g_n)) = (1,g_0,g_1,\dots,g_n)
$$
and is a p.s\. module homomorphism for each $n\ge 0$. The groups
$PH^*(G;M)\overset def\to{=} PH^*_G(EG_*;M)$ are
the p-bounded group cohomology groups of $G$ with coefficients in a p.s. $G$-module $M$.
The following result extends the Comparison Theorem of [O1]. When there is no
confusion, we write $EG_*$ for $(EG)_*$. However, for each $n$, $EG_n$ is the set $G^{n+1}$
while $(EG)_n = C_n(EG;\Bbb Q)$ is the free $\Bbb Q$-module on $EG_n$.

\proclaim{\bf\underbar{Theorem 1.1.10 - Comparison Theorem}} Let $(R_*,d_*)$ be a p.f\.
resolution of $\Bbb Q$ over $\Bbb Q[G]$ and $M$ a p.s\. $G$-module. Then there is an
isomorphism
$$
PH^*_G(EG_*;M)\cong PH^*_G(R_*;M) = H^*(PHom_G(R_*,M),\delta^*)
$$
\endproclaim

\prf As in [O1], one forms the bi-complex $EG_*\otimes R_*$.
Write $(EG)_p$ as $\Bbb Q[G][S_p]$ and $R_q$ as $\Bbb Q[G][T_q]$. The method of proof
of Proposition 1.1.6 provides isomorphisms
$$
PHom(\Bbb Q[S_p]\otimes \Bbb Q[G][T_q],M)
\cong PHom_G(\Bbb Q[G][S_p]\otimes \Bbb Q[G][T_q],M)
\cong PHom(\Bbb Q[G][S_p]\otimes \Bbb Q[T_q],M)
\tag1.1.11
$$
where the $G$-action on $\Bbb Q[G][S_p]\otimes \Bbb Q[G][T_q]$ is the diagonal one. Now
consider the bicocomplex formed by applying $PHom_G(\,\,\,,M)$ to $EG_*\otimes R_*$.
The $q^{th}$ row is
$$
PHom_G(EG_*\otimes R_q,M) = PHom_G(EG_*\otimes \Bbb Q[G][T_q],M)
$$
where the differential is the identity on the second coordinate. By the second isomorphism in
(1.1.11), the cohomology of this cocomplex is equal to the cohomology of
the cocomplex
$$
PHom(EG_*\otimes \Bbb Q[T_q],M)
$$
The standard p-bounded chain contraction on $EG_*$ yields a cocontraction of this cocomplex above
dimension zero. The resulting cohomology groups are zero in positive dimensions, with
$$
PH^0_G(EG_*\otimes \Bbb Q[T_q],M) = PHom(\Bbb Q[T_q],M)
\cong PHom_G(\Bbb Q[G][T_q],M) = PHom_G(R_q,M)
$$
Hence filtration by rows produces a spectral sequence with $E_1^{0,*} = PHom_G(R_*,M)$,
$E_1^{p,*} = 0$ for $p > 0$, and so $E_2^{0,*} = PH^*_G(R_*;M)$,
$E_2^{p,*} = 0$ for $p > 0$. Filtering by columns instead of rows reverses the roles of
$EG_*$ and $R_*$, resulting in a spectral sequence with $E_2^{*,0} = PH^*_G(EG_*;M)$,
$E_2^{*,q} = 0$ for $q > 0$.\hfill //
\endpf

The next result will provide our main technical tool for studying
the p-bounded group cohomology of a group $G$ with coefficients in a p.s\. $G$-module.

\proclaim{\bf\underbar{Theorem 1.1.12 - Serre Spectral Sequence}} Let
$(K,L_K)\inj (G,L_G)\surj (N,L_N)$ be a short-exact sequence of groups with word-length and
$M'$ a p.s\. $N$-module ($M'$ is then also a p.s\. $G$-module via the
surjection $G\surj N$). In addition, $M'$ is required to satisfy the hypothesis (1.1.H)
stated below. Then there is a first quadrant spectral sequence with
$$
E_2^{*,*} = \{PH^p_{N}(EN_*;PH^q(BK_*;M'))\}_{p,q\ge 0}
$$
converging to $PH^*_G(EG_*;M')$, with the natural transformation $PH^*(\,\,\,)\to
H^*(\,\,\,)$ inducing a map of Serre spectral sequences in cohomology.
\endproclaim

\prf As above $(EG)_n = \Bbb Q[G]^{\otimes n+1}$ with the $\Bbb Q[G]$-module structure induced by left
multiplication by $G$ on the left-most coordinate. Tensoring over $\Bbb Q[K]$ with $\Bbb Q$
yields the complex
$$
\Bbb Q[N]\leftarrow \Bbb Q[N]\otimes \Bbb Q[G]\leftarrow \Bbb Q[N]\otimes \Bbb Q[G]^{\otimes 2}
\leftarrow \dots
\tag1.1.13
$$
By the Comparison Theorem above, there are isomorphisms
$$
PH^*(BK_*;M')\cong PH^*_{K}((EK)_*;M')\cong PH^*_{K}((EG)_*;M')\cong PH^*(K\backslash(EG)_*;M')
$$
where $BK_* = K\backslash (EK)_*$ and $K\backslash(EG)_*$ is the complex in
(1.1.13).
Form the bicomplex
$B_{*,*} = EN_*\otimes K\backslash(EG)_*$. We will abbreviate the $(p,q)^{th}$ term of this
bicomplex as $N_p\otimes M_q$
where $N_p = \Bbb Q[N]^{\otimes p+1}$, $M_q = \Bbb Q[N]\otimes \Bbb Q[G]^{\otimes q}$.
Applying $PHom_N(\,\,\,,M')$ and filtering by rows produces a spectral sequence
which, by the isomorphisms of (1.1.11), collapses at the
$E_1^{*,*}$-term, with the only non-zero groups being
$$
E_1^{0,*} = PHom^*_N(K\backslash(EG)_*,M') = PHom^*_G((EG)_*,M')
$$
Computing $E_2^{*,*}$ yields
$E_2^{0,*} = PH^*_G(EG_*;M'), E_2^{p,*} = 0$ for $* > 0$.
We now consider the spectral sequence arising from filtration by columns. To compute the
$E_1^{*,*}$-term, we observe that the $p^{th}$ column is the cocomplex
$$
\dots\to PHom_N(N_p\otimes M_{q-1},M')\overset 1\otimes\delta^{q-1}\to
\longrightarrow PHom_N(N_p\otimes M_{q},M')\overset 1\otimes\delta^{q}\to
\longrightarrow PHom_N(N_p\otimes M_{q+1},M')\to\dots
\tag1.1.14
$$
The
$\Bbb Q[N]$-module structure on $N_p$ and $M_q$ is given by left multiplication by $N$ in the
left-most coordinate, and the $\Bbb Q[N]$-module structure on the tensor product is the
diagonal one. In order to properly identify the cohomology of the sequence in (1.1.14), we will
want to take partial adjoints.

\proclaim{\bf\underbar{Lemma 1.1.15}} For p.f\. $N$-modules $P$, $P'$, and p.s\. $N$-module $M'$
there are natural isomorphisms of p.s\. $G$-modules resp. p.s\. modules
$$
\gather
PHom(P\otimes P',M')\cong PHom(P,PHom(P',M'))\\
PHom_N(P\otimes P',M')\cong PHom_N(P,PHom(P',M'))
\endgather
$$
\endproclaim

\prf Again, these isomorphisms are well-known in the non p-bounded case.
We write $T$ resp. $T'$ for the orbit of $S$ resp. $S'$ under $G$.
As vector spaces, $P = \Bbb Q[T]$, $P' = \Bbb Q[T']$. Then
$P\otimes P' = \Bbb Q[T\times T']$ with weight function determined by setting
$|(t,t')| = |t| + |t'|$. The $G$ action on $T\times T'$ is the diagonal one.
For an element $f\in PHom(P\otimes P',M')$, denote its partial adjoint on the right by $\wt f$.
Thus $\wt f(gs)(g's') = f(gs,g's')$. Now suppose $\wt f : P\to PHom(P',M')$ is p-bounded.
Then there exists $C_1,n > 0$ depending only on $\wt f$ and $C_2$ depending only on $(x',U')$ such that
$$
\eta_{g(x',U')}(\wt f)(a))\le C_1C_2(1+|a|)^n(1+L(g))^n
$$
for all $g\in G$ and $a\in P$. Given $(t,t')\in T\times T'$, set $a=t$ and $U' = \{gt'\}$. Then
$$
\eta_{gx'}(f(t,t')) = \eta_{g(x',\{gt'\})}(\wt f)(t))\le C_1C_2(1+|t|)^n(1+L(g))^n
$$
implies $f$ is p-bounded on $T\times T'$, hence p-bounded. In the other direction,
suppose $f$ is p-bounded. As before, there are $C'_1,n' > 0$ depending only on $f$ and $C'_2$ depending only on $x'$ such that
$$
\eta_{gx'}(f(b))\le C'_1C'_2(1 + |b|)^{n'}(1+L(g))^{n'}
$$
for all $g\in G$ and $b\in P\otimes P'$.  Then
$$
\gather
\eta_{g(x',U')}(\wt f(t))\\
= \sum_{t_i'\in U'} \eta_{gx'}(f(t,g^{-1}t_i'))\\
\le \sum_{t_i'\in U'}  C_1'C_2'(1 + |(t,g^{-1}t_i')|)^{n'}(1+L(g))^{n'}\\
\le D_1D_2(1 +|t|)^n(1 + L(g))^{2n'}
\endgather
$$
where $D_1 = C_1'(\sum_{t_i'\in U'}(1+|t_i'|))^{n'}$ is independent of $g$ and $x'$ and $D_2 = C'_2$. This implies $f$ is p-bounded, which verifies the first isomorphism. The second follows from the first by the fact that the two adjoint maps preserve the $G$-action, hence induce isomorphisms on fixed-point sets.\hfill //
\endpf

Accordingly we may rewrite (1.1.14) as
$$
\aligned
\dots\to PHom_N(N_p, PHom(M_{q-1},M'))&\overset (\delta^{q-1})^*\to
\longrightarrow PHom_N(N_p, PHom(M_{q},M'))\\
&\overset (\delta^{q})^*\to
\longrightarrow PHom_N(N_p, PHom(M_{q+1},M'))\overset (\delta^{q+1})^*\to\longrightarrow
\dots
\endaligned
\tag1.1.16
$$
where $(\delta^k)^*$ is the map induced by $\delta^k : PHom(M_{k},M')\to PHom(M_{k+1},M')$.
Both $im(\delta^{k-1})$ and $ker(\delta^k)$ are
submodules of $PHom(M_{k},M')$ closed under the action of $N$, and so inherit a p.s\. $N$-module structure via
restriction (the p.s\. $N$-module structure on $PHom(M_{k},M')$ is that
given by Proposition 1.1.5).

\proclaim{\bf\underbar{Proposition 1.1.17}} There are equalities of p.s\.
$N$-modules
$$
\gathered
ker((\delta^k)^*) = PHom(N_p, ker(\delta^k))\\
im((\delta^{k-1})^*) = PHom(N_p, im(\delta^{k-1}))
\endgathered
\tag1.1.18
$$
\endproclaim

\prf First, $f : N_p\to PHom(M_k,M')$ maps to zero under $(\delta^k)^*$
exactly when $im(f)$ lies in $ker(\delta^k)$. Secondly,  $f\in
im((\delta^{k-1})^*)$ iff
there exists $f'\in PHom(N_p,PHom(M_{k-1},M'))$ with $f = f'\circ
\delta^{k-1}$.
But $f'\circ \delta^{k-1}$ is a map from $N_p$ to $im(\delta^{k-1})$. This
verifies the two equalities on the level of $N$-modules. They are
equalities of p.s\. $N$-modules because the p.s\. $N$-module structure on
both sides is induced by the restriction of a single p.s\. $N$-module
structure on $PHom(N_p, PHom(M_{k},M'))$.\hfill //
\endpf

In order to have an identifiable $E_1^{**}$-term, we need an additional hypothesis.
For applications below, we will state it in terms of a collection of hypotheses
indexed on the non-negative integers.

\underbar{Hypothesis 1.1.H(k)} For fixed $k\ge 0$,
$ker(\delta^k)/im(\delta^{k-1}) = PH^k(K;M')$ admits a p.s\. $N$-module structure for which
$$
ker(\delta^k)\surj PH^k(K;M')
$$
is an admissible epimorphism.
\vskip.05in

\underbar{Hypothesis 1.1.H} Hypothesis 1.1.H(k) is true for all $k\ge 0$.
\vskip.1in

Given that $im(\delta^{k-1})\inj ker(\delta^k)$ is an admissible monomorphism with the
p.s\. $N$-module structures as given above, this hypothesis is equivalent to the
statement that
$$
im(\delta^{k-1})\inj ker(\delta^k)\surj PH^k(K;M')
$$
is a short-exact sequence of p.s\. $N$-modules.
Under these conditions, Lemma 1.1.8 implies there is a corresponding short-exact sequence
$$
PHom_N(N_p, im(\delta^{k-1}))\inj PHom_N(N_p, ker(\delta^k))
\surj PHom_N(N_p, PH^k(K;M'))
$$
which together with (1.1.18) imply the $E_1^{*,*}$ is given by
$$
E_1^{p,q} = PHom_N(N_p, PH^q(K;M'))
$$
The $E_2^{**}$-term indicated in the statement of the theorem then follows as in the standard
Serre spectral sequence.\hfill //
\enddemo

A discrete group with word-length function $G$ has \underbar{p-bounded $A$-cohomology}
(where $A$ is a p.s\. $G$-module) if the natural transformation of cohomology
theories $PH_G^*(EG_*;A)\to H_G^*(EG_*;A)$ is an isomorphism. It is natural to suppose that the
class of groups with p-bounded cohomology is closed under arbitrary extensions. The following
corollary gives a partial result in this direction.

\proclaim{\bf\underbar{Corollary 1.1.19}} Let $K\inj G\surj N$ be a short-exact sequence of
groups equipped with word-length function, and let $A$ be a p.s\. $N$-module, such that hypothesis (1.1.H) is satisfied. If $K$ has p-bounded
cohomology with coefficients in $A$, and $N$ has p-bounded cohomology with coefficients in $PH^i(K;A) = H^i(K;A)$ for all $i$, then $G$ has p-bounded cohomology with coefficients in $A$.
\endproclaim

\prf The natural transformation from p-bounded cohomology to cohomology with coefficients in $A$
induces a map of Leray-Serre spectral sequences. With the given hypothesis,
there is an isomorphism of $E_2^{*,*}$-terms, where the $E_2^{**}$-term for p-bounded
cohomology is given in Theorem 1.1.12. By spectral sequence comparison, the result follows.
\hfill //
\endpf

Before giving the main application of this spectral sequence, we will
need a technical lemma.

\proclaim{\bf\underbar{Lemma 1.1.20}} Hypothesis 1.1.H(k) is satisfied
for $k = 0$ and $k = 1$.
\endproclaim

\prf When $k = 0$, $\delta^{k-1} = \delta^{-1} = 0$ and
$ker(\delta^k)\surj PH^k(K;M')$ is an isomorphism. So Hypothesis
1.1.H(0) is trivially satisfied. To handle the case $k = 1$, we
recall that the inclusion $K\inj G$ induces a p-bounded inclusion
of complexes
$$
(BK)_* = K\backslash (EK)_*\hookrightarrow K\backslash (EG)_*
\tag1.1.21
$$

\proclaim{\bf\underbar{Claim 1.1.22}} The inclusion of (1.1.21) induces an admissible
epimorphism of p.s\. cocomplexes
$$
(M^*,\delta^*) = (PHom(K\backslash (EG)_*,M'),\delta^*)
\surj (PHom(K\backslash (EK)_*,M'),\delta_K^*)
\tag1.1.23
$$
\endproclaim

\prf Let $\iota : S\hookrightarrow S'$ be an inclusion of sets. We also assume given
$\Bbb R^+$-valued maps $f_S$, $f_{S'}$ with $f_S = f_{S'}\circ \iota$.
Let $P$ resp. $P'$
be the p.f\. module generated by $(S,f_S)$ resp. $(S',f_{S'})$. We also suppose given a p-bounded
surjection $p : S'\surj S$ with $p\circ\iota = id$.
This surjection induces an admissible epimorphism $P'\surj P$ which we also denote by $p$.
Then for any p.s\. module $M'$, $\iota$
induces an admissible epimorphism $\iota^* : PHom(P',M')\surj PHom(P,M')$ with section
equal to $p^*$. In fact, if $\phi : P\to M'$ is p-bounded, then so is
$\phi' = \phi\circ p : P'\to M'$, and $\phi'\in PHom(P',M')$ maps to $\phi$ under $\iota^*$,
proving the surjectivity of $\iota^*$. For $\alpha\in PHom(P',M')$ the equality
$\eta_{(x,U)}(\iota^*(\alpha)) = \eta_{(x,\iota(U))}(\alpha)$
implies the p-boundedness of $\iota^*$. In the other direction the sequence
$$
\eta_{(x,U')}(p^*(\beta))
= \sum_{s_i'\in U'}\eta_x(p^*(\beta)(s_i'))
= \sum_{s_i'\in U'}\eta_x(\beta(p(s_i')))
\le C_2\eta_{(x,U)}(\beta)
$$
where $U = p(U')$ and $C_2 = |U| < \infty$
implies the p-boundedness of $p^*$. Thus $\iota$ is an admissible epimorphism.

Returning to the short-exact sequence
$K\overset i\to\inj G\overset p\to\surj N$, we fix a bounded section of sets
$s : N\inj G$; given $g\in G$ we denote the product $g(s(p(g)))^{-1}$ by $\lambda(g)$. 
For each $q\ge 0$ we have an inclusion of sets
$$
i_q : EK_q\inj EG_q
$$
induced by $i$ and a projection of sets
$$
\gather
p_q : EG_q\surj EK_q\\
(g_0,g_1,\dots ,g_q)\mapsto
(\lambda(g_0),\lambda(g_0)^{-1}\lambda(g_0g_1),\dots,
\lambda(g_0g_1\dots ,g_{q-1})^{-1}\lambda(g_0g_1\dots g_q))
\endgather
$$
The following properties are easily verified:
\vskip.05in

i) $i_q$ and $p_q$ are equivariant w.r.t. left multiplication by $K$ in the left-most
coordinate, hence descend to maps
$$
\gather
\ov i_q : K\backslash EK_q = BK_q\hookrightarrow K\backslash EG_q\\
\ov p_q : K\backslash EG_q\surj BK_q
\endgather
$$
\vskip.05in

ii) $p_q\circ i_q = id$ for each $q$, implying $\ov p_q\circ \ov i_q = id$ for each $q$.
\vskip.05in

iii) $\ov p_* = \{\ov p_q\}$ and $\ov i_* = \{\ov i_q\}$ are chain maps.
\vskip.05in

Now $(BK)_q = \Bbb Q[BK_q]$ and $(K\backslash EG)_q = \Bbb Q[K\backslash
EG_q]$ are the p.f\. modules generated respectively by the weighted sets
$(K^q,L(K\backslash K)_q)$ and $(K\backslash G\times G^q,L(K\backslash
G)_q)$, where
$$
\gathered
L(K\backslash K)_q(k_1,\dots k_q) = 1 + \sum_{i=1}^q L_K(k_i)\\
L(K\backslash G)_q(\ov g_0,g_1,\dots,g_q) = 1 + \ov L_G(g_0) + \sum_{i=1}^q
L_G(g_i)\,,\qquad \ov L_G(g_0) = \underset {k\in K}\to{min}\{L_G(kg_0)\}
\endgathered
$$
By what we have shown above, we conclude that
for any p.s\. $N$-module $M'$, the map
$$
(\ov i_q)^* : PHom((K\backslash EG)_q,M')\surj PHom((BK)_q,M')
$$
is an admissible epimorphism of p.s\. modules for each $q$ with section given by $p_q^*$.
As both $\ov i_*$ and $\ov p_*$ are chain maps, their duals $(\ov i_*)^*$ and $(\ov p_*)^*$
are cochain maps, which proves the claim.\hfill //
\endpf

We consider the following diagram
$$
\diagram
ker(\delta_K^1)\rto<-.5ex>|<<\tip_{p(1)}\ddouble &
ker(\delta^1)\lto<-.5ex>|>>\tip_{i(1)}\dto|>>\tip\\
ker(\delta_K^1)/im(\delta_K^0)\rto<-.5ex>|<<\tip_{\ov p(1)} &
ker(\delta^1)/im(\delta^0)\lto<-.5ex>|>>\tip_{\ov i(1)}
\enddiagram
\tag1.1.24
$$
where $\delta_K^*$ resp. $\delta^*$ are the coboundary maps appearing in
(1.1.23), and the surjections resp. injections in the diagram are those
induced by $(\ov i_*)^*$ resp. $(\ov p_*)^*$. The epimorphism $(\ov i_*)^*$ induces an 
isomorphism in cohomology by the Comparison Theorem, and $(\ov i_*)^*\circ (\ov p_*)^* =
id$, implying $\ov i(1)$ and $\ov p(1)$ are isomorphisms, and inverses of each other.
Also
$\delta_K^0 = 0$ since the $K$-module structure on $M'$ is trivial,
implying the left vertical map is the identity as indicated.
Denote the composition $p(1)\circ i(1)$ by $Pr$. Then for $f\in \ker(\delta^1)$,
$Pr(f)$ is given by the formula
$$
Pr(f)(\ov{g_0},g_1) = f(1,\lambda(g_0)^{-1}\lambda(g_0g_1))
$$
where $\ov{g_0}$ denotes the equivalence class of $g_0$ in $K\backslash G = N$. Now
$PHom(\Bbb Q[N\times G],M')$ is a p.s\. $N$-module with semi-norms indexed on the set
$\Lambda_{M'}\times \Cal P(N\times G)$.
Both $ker(\delta^1)$
and $im(\delta^0)$ inherit a p.s\. $N$-module structure via the inclusion into this
p.s\. $N$-module, inducing a quotient p.s\. $N$-module structure on 
$ker(\delta^1)/im(\delta^0)$.
Denoting the equivalence class of $f$ in this quotient by $[f]$, we have
$$
\eta_{(x,U)}([f])
= min\{\eta_{(x,U)}(f')\ |\ [f'] = [f]\}
$$
From the commutativity of the above diagram, we see that $[f'] = [f]$ implies $Pr(f') = Pr(f)$.
We claim that the map $[f]\mapsto Pr(f)$ is a monomorphism of p.s\. modules. In fact,
we have an inequality
$$
\eta_{(x,U)}(Pr(f))\le |U|\eta_{(x,Pr(U))}(f) = |U|\eta_{(x,Pr(U))}([f])
$$
where $Pr(U)$ denotes the image of $U\subset N\times G$ under the
composition $N\times G\overset {\ov p}_1\to\surj \{id\}\times
K\inj N\times G$; this implies the result. Since we have already
shown it is a splitting of the surjection $ker(\delta^1)\surj
ker(\delta^1)/im(\delta^0)$, we conclude that this surjection is
an admissible epimorphism of p.s\. $N$-modules, completing the
proof of Lemma 1.1.20.\hfill //
\endpf

An immediate consequence of this Lemma is the following
$5$-term sequence in p-bounded group cohomology.

\proclaim{\bf\underbar{Theorem 1.1.25 \ $5$-term sequence}} Let
$(K,L_K)\inj (G,L_G)\surj (N,L_N)$
be a short-exact sequence of groups with word-length, and $M'$ a p.s\. $N$-module.
Then there is a short-exact sequence
$$
0\to PH^1(N;M')\to PH^1(G;M')\to PH^0(N;PH^1(K;M'))\to PH^2(N;M')\to PH^2(G;M')
\tag1.1.26
$$
\endproclaim

\prf The proof follows exactly as in ordinary group cohomology by Lemma
1.1.20 and Theorem 1.1.12.\hfill //
\endpf
\vskip.3in

%%%%%%%%%%%%%%%%%%%%%%%%%%%%%%%%%%%%%%%%%%%%%%%%%

\subhead 1.2 The obstruction to injectivity
\endsubhead
\vskip.2in

Recall that the natural transformation $PH^*(G;\Bbb Q)\to H^*(G;\Bbb Q)$ induces a map of
dual vector spaces $(H^*(G;\Bbb Q))^*\to (PH^*(G;\Bbb Q))^*$, and that a group satisfies
property (PC3) if the composition
$$
\alpha_*(G) : H_*(G;\Bbb Q)\to (H^*(G;\Bbb Q))^*\to (PH^*(G;\Bbb Q))^*
\tag1.2.1
$$
is injective for all $*\ge 0$.

In this section we will work with a short-exact sequence
of groups with word-length
$$
(F',L_{F'})\inj (F,L_F)\surj (G,L_G)
\tag1.2.2
$$
where $F= F_S$ is a free group on a generating set $S$, $L_F$ is the word-length metric induced
by a function $f_S : S\to \Bbb N^+$ and (1.2.2) is the short-exact sequence of groups with
word-length associated to the short-exact sequence of groups $F'\inj F\surj G$ and the
word-length metric $L_F$.
The
Serre spectral sequence in homology with coefficients in a module $M$ has $E^2_{*,*}$-term
$$
E^2_{p,q} = H_p(G;H_q(F';M))
$$
and converges to $H_*(F;M)$. In general, $M$ is a non-trivial $F$-module; we use the same
notation for the $E^2_{*,*}$-term whether or not the corresponding action of $G$ on
$H_q(F';M)$ is trivial. Now $H_*(F;M) = H_*(F';M) = 0$ for $* > 1$. From this vanishing
we conclude

\proclaim{\bf\underbar{Proposition 1.2.3}} For all coefficient modules $M$, there are
isomorphisms
$$
H_p(G;M)\overset\cong\to\longrightarrow H_{p-2}(G;H_1(F';M))\quad p\ge 3
$$
and an injection
$$
H_2(G;M)\inj H_0(G;H_1(F';M))
$$
\endproclaim

We define $G$-modules inductively as follows:

\vskip.05in
(1.2.4.i) $B_0 = \Bbb Q$ with trivial $G$-module structure,
\vskip.05in

(1.2.4.ii) $B_m = H_1(F';B_{m-1})$, with diagonal $G$-module structure.
\vskip.05in

Note that $B_1 = H_1(F';\Bbb Q)$, and for $m > 1$ there is a natural isomorphism
$$
B_m\cong H_m((F')^m;\Bbb Q)\cong \otimes^m H_1(F';\Bbb Q)
$$
equipped with the diagonal conjugation action of $F$.
The action of $F$ on $B_m$ induced by this action of $F$ on $\otimes^m H_1(F';\Bbb Q)$
factors by the projection to $G$.

\proclaim{\bf\underbar{Proposition 1.2.5}} There are isomorphisms
$$
H_{2m}(G;B_n)\overset\cong\to\longrightarrow H_{2m-2}(G;B_{n+1})\quad m\ge 2
$$
and an injection
$$
H_2(G;B_n)\inj H_0(G;B_{n+1})
$$
\endproclaim

\prf This is a direct application of the previous proposition with $M = B_n$, as
$B_{n+1} = H_1(F';B_n)$.
\hfill //
\endpf

Starting at $H_{2m}(G;\Bbb Q)$, this proposition produces a sequence
$$
H_{2m}(G;\Bbb Q)\cong H_{2m-2}(G;B_1)\cong\dots \cong H_2(G;B_{m-1})\inj H_0(G;B_m)
\tag1.2.6
$$
where the maps in the sequence arise as differentials in the $E^2_{*,*}$-term of the
appropriate Serre spectral sequence.

A similar result holds for cohomology.

\proclaim{\bf\underbar{Proposition 1.2.7}} For all coefficient modules $M$, there are
isomorphisms
$$
H^{p-2}(G;H^1(F';M))\overset\cong\to\longrightarrow H^p(G;M) \quad p\ge 3
$$
and a surjection
$$
H^0(G;H^1(F';M))\surj H^2(G;M)
$$
\endproclaim

Let $B_m^* = Hom_{\Bbb Q}(B_m,\Bbb Q)$ denote the dual of $B_m$, with $G$-module structure
given by $gh(x) = h(g^{-1}x)$. The dual of Proposition 1.2.5 is

\proclaim{\bf\underbar{Proposition 1.2.8}} There are isomorphisms
$$
H^{2m-2}(G;B_{n+1}^*)\overset\cong\to\longrightarrow H^{2m}(G;B_n^*) \quad m\ge 2
$$
and a surjection
$$
H^0(G;B_{n+1}^*)\surj H^2(G;B_n^*)
$$
\endproclaim

This yields a sequence
$$
H^{0}(G;B_m^*)\surj H^{2}(G;B_{m-1}^*)\cong H^4(G;B_{m-2}^*)\cong\dots\cong  H^{2m}(G;\Bbb Q)
\tag1.2.9
$$
Again, the maps in the sequence occur as differentials in the $E_2^{*,*}$-term of the
appropriate Serre spectral sequence.

Denote the composition in (1.2.6) by $i_m$ and the composition in (1.2.9) by $j_m$. The
following commuting diagram derives from standard properties of the Serre spectral sequence.

$$
\diagram
H_{2m}(G;\Bbb Q) \rto|<<\tip^{i_m}\dto & H_0(G;B_m)\dto\\
(H^{2m}(G;\Bbb Q))^*\rto|<<\tip^{(j_m)^*} & (H^0(G;B^*_m))^*
\enddiagram
\tag1.2.10
$$

Let $A^*_0 = \Bbb Q$, and inductively set
$$
A^*_m = PH^1(F';A^*_{m-1})
\tag1.2.11
$$
for $m\ge 1$. For each $m$, a p.s\. $G$-module structure on $A^*_{m-1}$ induces a p.s\. $G$-module structure on
$A^*_m$, as shown in Lemma 1.1.20. This gives it a p.s\. $F$-module structure via the projection
$F\surj G$, which when restricted to $F'$ produces a p.s\. $F'$-module structure where the
action of $F'$ on both the module and indexing set is trivial. We will examine this
structure in more detail later on,
noting for now only its existence. Thus starting with the trivial p.f\. $G$-module structure
on $A^*_0$ as indicated following (1.1.4), we get a p.s\. $G$-module structure on each $A^*_m$.

\proclaim{\bf\underbar{Lemma 1.2.12}} For all $m\ge 0$ there is a Serre spectral sequence
$$
E_2^{p,q} = PH^p(G;PH^q(F';A_m^*))
$$
converging to $PH^{p+q}(F;A_m^*)$. Moreover, $PH^q(F';A_m^*) = 0$ for $q > 1$.
\endproclaim

\prf
Recall $F$ is a free group with basis $S$ where it is assumed that $S\cap S^{-1}=\emptyset$. 
There is a $\Bbb Q[F]$-free resolution of
$\Bbb Q$
$$
M_S\overset\eta\to\to \Bbb Q[F]\to \Bbb Q
\tag1.2.13
$$
where $M_S = \Bbb Q[F][S]$ is the free $\Bbb Q[F]$-module on $S$, $\Bbb Q[F]$ the free 
$\Bbb Q[F]$-module on the single generator $[1]$, and $M_S\to \Bbb Q[F]$ is the
$\Bbb Q[F]$-module map defined on basis elements by $\eta([x]) = (x-1)$. The map
$\eta$ induces an isomorphism between $M_S$ and the augmentation ideal $I[F]$. We denote
this \lq\lq short complex\rq\rq in (1.2.13) by $R_*(F)$.
For $x\in S\coprod S^{-1}$, set
$$
\psi(x) =
\left\{\aligned
[x]\,\,\text{if}\, x\in S\\
-x[x^{-1}]\,\,\text{if}\, x^{-1}\in S
\endaligned
\right.
$$

Each $g\in F$ admits a unique reduced word representation $g = x_1x_2\dots x_n$ where
$x_i\in S\coprod S^{-1}$ for each $i$. Define $s_0 : \Bbb Q[F]\to M_S$ as the
$\Bbb Q$-vector space map given on basis elements by
$$
s_0(g[1]) = s_0(x_1x_2\dots x_n[1]) = \sum_{j=1}^n x_1x_2\dots x_{j-1}\psi(x_j)
$$
The fact that $L_F$ is a word-length metric implies $| x_1x_2\dots x_{j-1}\psi(x_j)|\le 2L(g)$
for each $j$, and also that $L(g)\ge n$. From this, one concludes that for each $m\in\Bbb N$
there is a sequence of inequalities
$$
|s_0(g)|\le n(1 + 2L(g))\le (1 + 2L(g))^2
$$
which in turn implies $s_0$ is p.s\. module homomorphism with respect to the p.f\. $F$-module
semi-norms on $\Bbb Q[F]$ resp. $\Bbb Q[F][S]$. Since $\eta$ is easily seen to be a p.s\.
$F$-module homomorphism, we conclude that the short complex described above is a p.f\.
resolution of $\Bbb Q$ over $\Bbb Q[F]$.

\proclaim{\bf\underbar{Proposition 1.2.14}} Let $F'$ be a subgroup of $F$ with induced metric. Then
$$
PH^*(F';A)\cong PH^*_{F'}(R_*(F);A)
$$
for any p.s\. $F'$-module $A$.
\endproclaim

\prf  As the function $f_S$ is $\Bbb N^+$-valued, so is the word-length metric $L_F$.
Let $T = F'\backslash F$ denote the right coset space, and $p : F\surj T$ the natural
projection. Define $f_T : T\to \Bbb R^+$ by $f_T(F'g) = min\{L_F(f'g)\,|\, f'\in F'\}$.
Note that as $L_F$
is $\Bbb N^+$-valued, so is $f_T$. Thus for all $g\in F$, there exists an $f'_g\in F'$ with
$f_T(F'g) = L_F(f'_g g)$. Choosing such an $f'_g$ for each $g$ and writing
$F'g\in T$ as $\ov g$, we set $s(\ov g) = f'_g g$. Also, we will write $f_{T\times S}$ for
the function $T\times S\ni (t,s)\mapsto f_T(t) + f_S(s)$. Finally
for $g\in F$ we denote $g(s(\ov g))^{-1}\in F'$ by $\lambda(g)$, as before. By construction,
$L_F(\lambda(g))\le 2L_F(g)$. We now consider the morphism of complexes
$$
\xymatrix{
\Bbb Q[F][S]\ar@<1ex>[r]^{d_1}\ar[d]^{\phi_1} & \Bbb Q[F]\ar[d]^{\phi_0}\ar@<1ex>[l]^{s_0}\\
\Bbb Q[F'][T\times S]\ar@<1ex>[r]^{d_1'} & \Bbb Q[F'][T]\ar@<1ex>[l]^{s_0'}}
\tag1.2.15
$$
The top row is $R_*(F)$. Denote the bottom row by $\wt R_*(F')$. The p.f\. $F'$-module
structure on $\Bbb Q[F'][T\times S]$ resp. $\Bbb Q[F'][T]$ is that induced by the function
$f_{T\times S}$ resp. $f_T$, and the left action of $F'$ on $\Bbb Q[F']$.
The maps $\phi_i$ and their inverses are defined as
$$
\gathered
\phi_1(g[s]) = \lambda(g)[\ov g,s],\qquad \phi_1^{-1}(g'[\ov g,s]) = g's(\ov g)[s]\\
\phi_0(g) = \lambda(g)[\ov g],\qquad \phi_0^{-1}(g'[\ov g]) = g's(\ov g)
\endgathered
\tag1.2.16
$$
It is easily verified that $\phi_i^{-1}$ is non-increasing in norm, while $\phi_i$ increases
the norm by no more than a factor of three. Defining
$$
\gather
d_1' = \phi_0\circ d_1\circ \phi_1^{-1}\\
s_0'=  \phi_1\circ s_0\circ \phi_0^{-1}
\endgather
$$
in diagram (1.2.15) makes $\phi_*$ a p-bounded $\Bbb Q[F']$-module isomorphism of complexes
with p-bounded inverse. Moreover, the p-boundedness of $\phi_*$ and $\phi_*^{-1}$ make
the contraction $s_0'$ p-bounded. Thus $\wt R_*(F')$ satisfies the hypothesis of the Comparison
Theorem. We then have isomorphisms
$$
PH^*(F';A)\cong PH^*_{F'}(\wt R_*(F');A)\cong PH^*_{F'}(R_*(F);A)
$$
by the Comparison Theorem, together with the p-boundedness of
$\phi_*$ and its inverse.\hfill //
\endpf

The complex $R_*(F)$ is zero above dimension one. Thus

\proclaim{\bf\underbar{Corollary 1.2.17}} For all (free) subgroups with word-length function
$(F',L_{F'})$ of $(F,L_F)$, ($L_{F'} = (L_F)|_{F'}$) and p.s. $F'$-modules $M$,
$PH^*(F';M) = 0$ for $*\ge 2$.
\endproclaim

We now return to the bicomplex used in the proof of Theorem 1.1.12. Referring to (1.1.18), we
see that Corollary 1.2.17 implies $ker(\delta^k) = im(\delta^{k-1})$ for all $k\ge 2$, so that
Hypothesis 1.1.H(k) is trivially satisfied for $k\ge 2$. By Lemma 1.1.20, Hypothesis 1.1.H(k)
is always satisfied for $k=0$ and $k=1$.
Consequently Hypothesis 1.1.H is satisfied as stated, and Theorem 1.1.12 applies,
completing the proof of Lemma 1.2.12.\hfill //
\endpf

\proclaim{\bf\underbar{Corollary 1.2.18}} There is a sequence
$$
PH^{0}(G;A_m^*)\surj PH^{2}(G;A_{m-1}^*)\cong PH^4(G;A_{m-2}^*)\cong\dots\cong  PH^{2m}(G;\Bbb Q)
\tag1.2.19
$$
where the maps in the sequence occur as differentials in the $E_2^{*,*}$-term of the
appropriate Serre spectral sequence for p-bounded cohomology.
\hfill\text{//}
\endproclaim

The proof is exactly as before, given the previous lemma. The natural transformation
$PH^*(G;A)\to H^*(G;A)$ induces an equally natural transformation
$$
(H^*(G;A))^*\to (PH^*(G;A))^*
$$
Together with (1.2.6)
and the duals of (1.2.9) and (1.2.19) we arrive at the following commuting diagram which is
an extension of (1.2.10)
$$
\diagram
H_{2m}(G;\Bbb Q)\rto|<<\tip \dto|<<\tip & H_0(G;B_m)\dto|<<\tip\\
(H^{2m}(G;\Bbb Q))^*\rto|<<\tip \dto & (H^0(G;B_m^*))^*\dto\\
(PH^{2m}(G;\Bbb Q))^*\rto|<<\tip & (PH^0(G;A_m^*))^*
\enddiagram
\tag1.2.20
$$
The injectivity of the horizontal arrows follows from what we have already shown.
Recall property (PC3) is the statement
$$
\alpha_*(G) : H_*(G;\Bbb Q)\to (H^*(G;\Bbb Q))^*\to (PH^*(G;\Bbb Q))^*
\tag1.2.21
$$
is injective for all $*\ge 0$.

\proclaim{\bf\underbar{Proposition 1.2.22}} If the composition
$$
\beta_m(G) : H_0(G;B_m)\inj (H^0(G;B_m^*))^*\to (PH^0(G;A_m^*))^*
$$
is injective for all $m\ge 0$, then $\alpha_{2m}(G)$ is injective for all $m\ge 0$.
\endproclaim

\prf This is an immediate consequence of (1.2.20).\hfill //
\endpf

%%%%%%%%%%%%%%%%%%%%%%%%%%%%%%%%%%%%%%%%%%%%%%%%%

\vskip.3in
\subhead
1.3 Analysis of the obstruction
\endsubhead
\vskip.2in

Proposition 1.2.22 above identifies a condition sufficient to guarantee injectivity of $\alpha_n(G)$ in even dimensions, and a similar analysis works in odd dimensions after crossing $G$ with $\Bbb Z$. 
The purpose of this section is to indicate the relationship between the injectivity of $\beta_m(G)$ and the first Dehn function of $G$ when $G$ is finitely-presented. In dimension 2 (cf. Theorem 2.1.3 below), $\alpha(G,L_G)$ is an isomorphism with arbitrary coefficients when the first Dehn function of $G$ is of polynomial type. In higher dimensions,  the injectivity of $\beta_m(G)$ follows if one can show that a certain natural class of projection maps are admissible epimorphisms (Theorem 1.3.5 below). This section is in preliminary form; a sequel to this paper will contain a much more detailed analysis of Conjecture A, along with complete proofs of the results stated in this section.
\vskip.2in

We begin by recalling the definition of Dehn functions. Let $\Cal P = <\Cal S\ |\ \Cal W>$ be a finite presentation of a discrete group $G$.
Then there is a short-exact sequence
$$
F'\inj F\surj G
\tag1.3.1
$$
where $F$ is the free group on the (finite) set of generators $\Cal S$ and $F'$ the
subgroup of $F$ normally generated by the (finite) set of relators $\Cal W$. We take
$L_F$ to be the standard word-length metric on $F$ which takes the value $1$ on each
generator, with $L_G$ the standard word-length function on $G$ induced by $L_F$ and the
projection $F\surj G$. For $w\in\Cal W$, denote its image in $F'$ by $\ov w$. Any element $y\in F'$ may be written as 
$$
y = {\ov w_1}^{x_1}{\ov w_2}^{x_2}\dots {\ov w_n}^{x_n}
\tag1.3.2
$$
where for each $i$, $w_i^{\pm 1} \in \Cal W$, $x_i\in F$ and ${\ov w}^x = x{\ov w}x^{-1}$. The \underbar{area} of $y$,
written $Area_{\Cal P}(y)$, is the minimum number $n$ such that $y$ can be written as in
(1.3.2). A map $f : \Bbb N\to \Bbb N$ is called an \underbar{isoperimetric function} for
the presentation if
$$
Area_{\Cal P}(y)\le f(n)
$$
for all relations $y$ with $L_F(y)\le n$. Among all isoperimetric functions associated to
$\Cal P$ there is a minimal one, $f_{\Cal P}$, referred to as the \underbar{Dehn function}
of the presentation $\Cal P$.

Dehn functions are due to Gersten ([Ger1], [Ger2]). We say that the Dehn function is of 
\underbar{polynomial type} if it is bounded above by a polynomial function.

Some notation. We will write
$H_1(F';\Bbb Q)$ as $R = R(1)$, and in general for $n \ge 1$ denote 
$\otimes^n H_1(F';\Bbb Q)$ by $R(n)$. Recall that as a subgroup of $F$, $F'$ is normally
generated by elements of the form $(\ov w)^g$ where $w\in \Cal W$ and $g\in F$. The image of $(\ov w)^g$ in $H_1(F';\Bbb Q)$ only depends on $\ov w$ and the image of $g\in G$. Thus $R$ is spanned as a vector space over $\Bbb Q$
by $\{[(\ov w)^g]\}_{w\in\Cal W,g\in G}$, where $[(\ov w)^g]$ denotes the image of
$(\ov w)^g\in F'$ under the canonical map $F'\to R$. From this we see there is a natural
surjection
$$
\Bbb Q[G][\Cal W]\surj R\,\,\,,\qquad (g,w)\mapsto [(\ov w)^g]
\tag1.3.3
$$
Let $P(n) =\otimes^n \left(\Bbb Q[G][\Cal W]\right)$. The map in (1.3.3) induces a surjection of $n$-fold tensor products
$$
p_n: P(n)\surj R(n)\,\,\,,\qquad 
((g_1,w_1),\dots,(g_n,w_n))\overset p_n\to\mapsto ([(\ov w_1)^{g_1}],\dots,[(\ov w_n)^{g_n}])
\tag1.3.4
$$
\vskip.2in
The word-length function on $F'$ induces a weight function on $\Cal W$; together with the 
word-length function on $G$ we get a p.f\. $G$-module structure on   $P(n)$. This induces 
a (quotient) p.s\. $G$-module structure on $R(n)$, where in both cases the $G$-action is 
the diagonal one. Let $T_n = \prod_1^n (G\times\Cal W)$, so that $P(n) = \Bbb Q[T_n]$. 
For each $x\in T_n$, let $[x]$ denote its image in $H_0(G;P(n)) = P(n)_G$, and $[p_n(x)]$ 
its image in $R(n)_G$. For $[p_n(x)]\ne 0$, let $\Bbb Z _{[p_n(x)]}$ be the copy of $\Bbb Z$ 
generated by this element. Again, the p.s\. $G$-module structure above induces a quotient 
p.s\. module structure on $R(n)_G$, and so by restriction a p.s\. module structure on 
$\Bbb Z _{[p_n(x)]}$ for each $[p_n(x)]\ne 0$. Finally for each such $[p_n(x)]$ we may 
restrict the p.s\. module structure on $R(n)_G$ to $\Bbb Z _{[p_n(x)]}$. It is not hard 
to show that this is the same as the p.f\. module structure induced by the (quotient) 
length function on $\Bbb Z _{[p_n(x)]}$.
Note that for each $x$  there is a canonical word-length metric on $\Bbb Z_{x} =$ the 
subgroup of $P(n)$ generated by $x$ and a projection $\Bbb Z_{x}\to \Bbb Z _{[p_n(x)]}$ 
which is an isomorphism of abelian groups
 
\proclaim{\bf\underbar{Theorem 1.3.5}} If the projection map 
$p_{x} :\Bbb Z_{x}\to\Bbb Z _{[p_n(x)]}$ is an admissible epimorphism for each $x$ with 
$[p_n(x)]\ne 0$, then the map $\beta_n(G)$ is injective.
\endproclaim

We give a sketch of the proof; a more detailed proof (including a proof of the two technical 
lemmas below) will appear in the sequel. The hypothesis of the theorem implies that the 
projection map $\Bbb Z_{x}\to\Bbb Z _{[p_n(x)]}$ is a p-equivalence (i.e., a p-bounded 
isomorphism of p.s\. modules). This in turn implies that the p-bounded homomorphism 
$\Bbb Z_{x}\to\Bbb Q$ induced by the inclusion of $\Bbb Z$ into $\Bbb Q$ (equipped with 
standard p.f\. module structure) factors by the projection $p_{x}$. Denote the p-bounded 
homomorphism $\Bbb Z _{[p_n(x)]}\to\Bbb Q$ by $\phi_{[p_n(x)]}$. Now let $V$ be a subspace 
of $R(n)_G$ spanned by a \underbar{finite} number of elements 
$\{[p_n(x_1)],[p_n(x_2)],\dots,[p_n(x_n)]\}$. Then $V$ inherits a p.s\. module structure 
from its embedding into $R(n)_G$. 

\proclaim{\bf\underbar{Lemma 1.3.6}} Let $0\ne a\in V$. Then there exists $1\le i\le n$ 
and a p-bounded extension $\phi:V\to\Bbb Q$ of 
$\phi_{[p_n(x_i)]}:\Bbb Z _{[p_n(x_i)]}\to\Bbb Q$  with $\phi(a)\ne 0$.
\endproclaim

\proclaim{\bf\underbar{Lemma 1.3.7}} If $V$ is a finite-dimensional subspace of $R(n)_G$ 
and $\phi :V \to\Bbb Q$ a p-bounded homomorphism, then $\phi$ extends to a p-bounded 
homomorphism $\phi : R(n)_G\to\Bbb Q$.
\endproclaim

The proof of Lemma 1.3.6 uses in an essential way the fact that $V$ is finite-dimensional, 
while Lemma 1.3.7 is an analogue of the Hahn-Banach theorem. Together they imply Theorem 1.3.5. For suppose $0\ne y\in H_0(G;B_n) = R(n)_G$. Because homology has finite supports, the image of $y$ lies in a finite-dimensional subspace of the type considered in Lemma 1.3.6, which guarantees the existence of a p-bounded homomorphism $\phi:V\to\Bbb Q$ with $\phi(y)\ne 0$. By Lemma 1.3.7, this homomorphism may be extended over all of $R(n)_G$, yielding an element of $PH^0(G;A_n^*)$ which pairs non-trivially with the image of $y$. This implies $\beta_n(G)(y)\ne 0$. Varying $y$ then implies the injectivity of $\beta_n(G)$.

One may consider a hypothesis similar to that in Theorem 1.3.5, but without passing to 
$G$-invariants. For $x\in T_n$, $p_n(x)$ denotes the image in $R(n)$. If $p_n(x)\ne 0$,
we denote by $\Bbb Z_{p_n(x)}$ the copy of $\Bbb Z$ in $R(n)$ generated by $p_n(x)$,
and by $\Bbb Z_{x}$ the corresponding copy of $\Bbb Z$ in $P(n)$ generated by $x$. 

\proclaim{\bf\underbar{Theorem 1.3.8}} Suppose that $G$ is a finitely-presented group with
Dehn function of polynomial type. Then for each $n\ge 1$ and $x\in T_n$ with $p_n(x)\ne 0$,
the projection map $p_n$ induces an equivalence of p.s\. modules 
$\Bbb Z_{x}\to \Bbb Z_{p_n(x)}$.
\endproclaim

In fact, the condition on the Dehn function of $G$ directly implies the result for $n=1$, and the result for $n > 1$ follows directly from the case $n=1$.

It follows from Theorem B of the introduction (proved below in section 2) that $\beta_1(G)$ is an injection when the first Dehn function is of polynomial type. In fact, using the five-term exact sequence the hypothesis of Theorem 1.3.5 is easily verified in this case. However, unlike the situation in Theorem 1.3.8, the case  $n > 1$ in Theorem 1.3.5 does not follow in any obvious way from the case $n = 1$.

%%%%%%%%%%%%%%%%%%%%%%%%%%%%%%%%%%%%%%%%%%%%%%%%%

\vskip.3in
\subhead
1.4 The map $\alpha_1(G,L_G)$
\endsubhead
\vskip.2in

As we have observed in the introduction, $\alpha_i(G, L_G)$ depends only on $G$ and the choice
of length function. Thus if $G$ is generated by a set $S$ and $L_G$ is the word-length function
determined by a function $f : S\to \Bbb N^+$ (i.e., by the word-length metric $L_F$ determined
by $f$ on the free group $F$ generated by $S$, together with the natural surjection 
$F\surj G$) then $L_G$ and so also $\alpha_i(G, L_G)$
is independent of the choice of relator set $W$ in a presentation $P = <S|W>$ of $G$ which 
uses $S$ as the set of generators. We verify Theorem B of the introduction by proving the 
following three Lemmas. We assume that $G_{ab}$, the abelianization of $G$, is finitely
generated.
\vskip.1in

\proclaim{\bf\underbar{Lemma 1.4.1}} Let $\ov G = G_{ab}/G_{ab}^{torsion}$. Then there is a
system of generators $\ov S$ for $\ov G$ and a weight function $\ov f : \ov S\to \Bbb N^+$
for which $\alpha_1(\ov G, L_{\ov G})$ is injective, where $L_{\ov G}$ is the word-length
function determined by $\ov f$.
\endproclaim

\prf By the assumption on $G_{ab}$, $\ov G\cong \Bbb Z^r$ for some finite integer $r \ge 0$.
Let $\ov S = \{\ov x_1,\dots \ov x_r\}$ be a basis for $\ov G$, and set
$f_{\ov S}(\ov x_i) = 1$ for all $i$.
Let $L_{\ov G}$ be the word-length function determined by $f_{\ov S}$ (in other words, 
the standard word-length function associated with this set of generators). Then
$\alpha_1(\ov G, L_{\ov G})$ is an isomorphism as observed above.\hfill //
\endpf

\proclaim{\bf\underbar{Lemma 1.4.2}} Given $\ov S$ and $f_{\ov S}$ as in the previous Lemma, 
there exists a generating set $S$ for $G$ and proper function $f : S\to \Bbb N^+$ so that
$(G,L_G)\surj (\ov G, L_{\ov G})$ is surjection of groups with word-length function.
\endproclaim

\prf We first choose a set of elements $S'\subset G$ which maps isomorphically to
$\ov S\subset \ov G$ under the surjection $G\surj \ov G$, and set
$f_{S'}(x'_i) = f_{\ov S}(\ov x_i)$ where
$x'_i\in S'$ maps to $\ov x_i$. Let $S''$ be an arbitrary set of generators for
$ker(G\surj \ov G)$ and let $f_{S''} : S''\to \Bbb N^+$ a proper function on $S''$. Then
$S = S'\coprod S''$ is a generating set for $G$ equipped with proper function
$f = f_{S'}\coprod f_{S''}$. Setting $L_G$ to be the word-length function determined by $f$ 
completes the proof.\hfill //
\endpf

\proclaim{\bf\underbar{Lemma 1.4.3}} Suppose $\phi :(G,L_G)\to(H,L_H)$ is a
p-bounded homomorphism of
groups with word-length function. If $\phi_1 : H_1(G;\Bbb Q)\to H_1(H;\Bbb Q)$ and
$\alpha_1(H,L_H)$ are injective, then $\alpha_1(G,L_G)$ is injective.
\endproclaim

\prf This follows from the naturality  of $\alpha_1(G,L_G)$ with
respect to p-bounded homomorphisms of groups equipped with
word-length function.\hfill //
\endpf

Taking $(H,L_H) = (\ov G,L_{\ov G})$, these three Lemmas together imply Theorem B.

\newpage
\vskip.2in

%%%%%%%%%%%%%%%%%%%%%%%%%%%%%%%%%%%%%%%%%%%%%%%%%
%%%%%%%%%%%%%%%%%%%%%%%%%%%%%%%%%%%%%%%%%%%%%%%%%

\centerline{\bf\underbar{Higher Dehn functions}}
\vskip.3in

\subhead
2.1 Dehn Functions and simplicial resolutions
\endsubhead
\vskip.2in

We begin by considering a variant of the Dehn function associated
to a presentation. As already noted, an element $w\in F'$ can be
written as
$$
w = w_1^{x_1}w_2^{x_2}\dots w_n^{x_n}
\tag2.1.1
$$
where for each $i$, $w_i^{\pm 1} \in \Cal W$ and $x_i\in F$, and
$Area_{\Cal P}(w)$ is the minimum number $n$ such that $w$ can be written as in
(2.1.1).
Analogously, define $Area'_{\Cal P}(w)$ as the smallest integer $m'$ such that
$w$ can be written as in (2.1.1) with $m' = \sum_{i=1}^k L(w_i) +
2L(x_i)$. Let $f'_{\Cal P}$ be the minimal isoperimetric function
defined using $Area'_{\Cal P}$ instead of $Area_{\Cal P}$. The
inequalities $Area_{\Cal P}\le Area'_{\Cal P}$ and $f_{\Cal P}\le
f'_{\Cal P}$ are obvious. The following result is due to Gersten
[GC].

\proclaim{\bf\underbar{Lemma 2.1.2}} Let $M = max\{L(w)\ |\ w\in \Cal W\}$. Then
$$
f'_{\Cal P}(n)\le (2M)f_{\Cal P}(n)^2 + (2n + M)f_{\Cal P}(n)
$$
\endproclaim

In preparation for what follows, we will need to recall some terminology
and constructions used in [O1]. For standard properties of simplicial sets, we refer
the reader to [M].
\vskip.1in

For a simplicial group $\Gamma\hskip-.03in .$ set
$$
\gathered
\Gamma^{-1}_n = \Gamma_n\\
\Gamma^k_n = \bigcap_{i=0}^k ker(\partial_i:\Gamma_n\surj \Gamma_{n-1})\quad\text{for}\,k\ge 0, n\ge 1\\
\Gamma^0_0 = \partial_1(\Gamma^0_1)
\endgathered
$$
For $0\le k < n$ and $n\ge 1$ there is a split short-exact sequence
$$
\Gamma^{k}_n\inj \Gamma^{k-1}_n\overset\partial_k\to\surj \Gamma^{k-1}_{n-1}
$$
with the splitting induced by the restriction of $s_k$ to $\Gamma^{k-1}_{n-1}$. When
$k = n$ there is an exact sequence
$$
\Gamma^{n}_n\inj \Gamma^{n-1}_n\overset\partial_n\to\to \Gamma^{n-1}_{n-1}
$$
and the Kan extension property implies
$$
\pi_n(\Gamma\hskip-.03in .) = \Gamma^n_n/(\partial_{n+1}(\Gamma^n_{n+1}))
$$
where $\Gamma\hskip-.03in .$ is viewed here as a simplicial set with basepoint
$1\in \Gamma_0$. One also has
$$
\pi_0(\Gamma\hskip-.03in .) = \Gamma_0/\Gamma^0_0
$$

We say that $\Gamma\hskip-.03in .$ is a \underbar{resolution} if
$\pi_n(\Gamma\hskip-.03in .) = 0$ for all $n > 0$.
This is equivalent to the condition that
$\Gamma^{n-1}_n\overset\partial_n\to\to \Gamma^{n-1}_{n-1}$ is a surjection
for all $n\ge 1$ (note that it need not be a split-surjection).

A simplicial group $\gm\hskip-.03in .$ equipped with word-length function $L.$ is a
simplicial group $\gm\hskip-.03in . = \{[n]\mapsto \gm_n\}_{n\ge 0}$ where $L_n$ is a
word-length
function on $\gm_n$ for each $n$, and all face and degeneracy maps are
\underbar{p-bounded}.
The simplicial group together with its word-length function will be
written as a pair $(\gm\hskip-.03in ., L.)$.
A word-length function on an augmented simplicial group
$\gm\hskip-.03in .^+$ is a word-length function $L.$ on the associated simplicial group
$\gm\hskip-.03in .$ ($\gm_n = \gm^+_n$ for $n\ge 0$)
together with a word-length function $L_{-1}$ on $\gm_{-1}$ induced by $L_0$ and the
augmentation map
$\varepsilon : \gm_0\surj \gm_{-1}$. The
resulting augmented simplicial group together with word-length function is
written as $(\gm\hskip-.03in .^+,L.^+)$. The associated simplicial group with
word-length function $(\gm\hskip-.03in .,L.)$ is gotten by restricting to simplicial
dimensions $n\ge 0$.
A map $\phi$ of simplicial or augmented simplicial groups is
\underbar{p-bounded} if it is p-bounded in each simplicial degree.
\vskip.1in

Occasionally we need to keep track of
generating sets, in which case they are included in the notation. As always, we assume
generating sets are countable. We call
$(F,\Bbb X,L)$ a \underbar{triple} when $F$ is the free group with basis $\Bbb X$
equipped with a function
$f : \Bbb X\to \Bbb N^+$, and $L_F$ is the word-length metric
induced by $f$. This definition extends to the augmented simplicial
setting. A triple $(\gm\hskip-.03in .^+,\Bbb X.^+,L.^+)$ indicates
i) $(\gm\hskip-.03in .^+, L.^+)$ is an augmented simplicial group with word-length,
ii) $(\gm_n,\Bbb
X_n,L_n)$ is a triple for each $n\ge 0$ and iii) $\gm_{-1}$ is generated by
$\Bbb X_{-1} = \Bbb
X_0$.
We do not put any additional restriction on the face and degeneracy
maps when including a generating set $\Bbb X.^+$ (although in practice it
can always be arranged for $\Bbb X.^+$ to be closed under degeneracies).

For an augmented simplicial group $\Gamma\hskip-.03in .^+$, we will denote the kernel
$\ker(\gm\hskip-.03in .^+\surj \gm_{-1})$ of the simplicial augmentation map as
$\gm(\varepsilon)\hskip-.015in .^+$. This is an augmented simplicial subgroup of
$\gm\hskip-.03in .^+$ with $\gm(\varepsilon)_{-1} = \{1\}$.

We say that a free resolution $(\gm\hskip-.03in .^+,L.^+)$ or $(\gm\hskip-.03in .^+,\Bbb
X.^+,L.^+)$ is \underbar{type $P(m)$} if $\gm(\varepsilon)\hskip-.015in .^+$, viewed as
a \underbar{simplicial set}, admits a simplicial contraction through dimension
$(m-1)$ which is p-bounded in each degree (with respect to $L.^+$).
The resolution is \underbar{type $P$} if $\gm(\varepsilon)\hskip-.015in .^+$ admits a
simplicial contraction (of simplicial sets) which is p-bounded in all
degrees. Type $P$ is slightly stronger than being type $P(m)$ for all $m$.

\proclaim{\bf\underbar{Theorem 2.1.3}} If $G$ is a finitely-presented group with polynomial
Dehn function $f_G$, then $\alpha_2(G,L^{st}_G)$ is an isomorphism.
\endproclaim

\prf
As above, we denote the finite set of generators of
$G$ by $\Cal S$ and the finite set of relators by $\Cal W$. Let $\gm_0$ be the free group
on $\Bbb X_0 = \Cal S$, and $\gm_1$ the free group on $\Bbb X_1 = \Cal S\coprod \Cal W$.
Let $\varepsilon : \gm_0\surj G$ be the obvious projection. The natural
inclusion
$\Bbb X_0\hookrightarrow \Bbb X_1$ determines a monomorphism
$s_0 : \gm_0\inj \gm_1$. Define
$\partial_i$ on $\gm_1$ ($i = 0,1$) as the unique homomorphism determined on generators
by
$$
\gathered
\partial_i(s) = s\,\,\text{if}\,s\in S, i = 0,1\\
\partial_0(w) = id, \quad \partial_1(w) = \ov w
\endgathered
\tag2.1.4
$$
where $\ov w$ denotes $w$ viewed as an element of $\gm_0$. For $m \ge 2$ let $\gm_m$
be the free group on
$$
\Bbb X_m = \coprod\{s(\Bbb X_1)\}
$$
where the coproduct is over all iterated degeneracies from dimension $1$ to dimension $m$.
Finally let $\gm_{-1} = G$. The partial simplicial structure on
$\{\gm_n\}_{-1\le n\le 1}$ defined above admits a unique extension
to an augmented simplicial structure on $\gm\hskip-.03in.^+ = \{\gm_n\}_{n\ge -1}$. The
word-length
function $L.^+ = \{L_n\}_{n\ge -1}$ on $\gm\hskip-.03in.^+$ is the standard one in
dimensions $-1$ and $0$.
In dimension $1$ it is the metric determined by the function
$x\mapsto L_0(x)$, $w\mapsto L_0(\ov w)$ where $x\in \Cal S$ and $w\in \Cal W$.
In dimensions $\ge 2$ it is the unique metric defined on generators by
$L_m(s(x))= L_1(x)$ where $x\in \Bbb X_1$ and $s$ is an iterated degeneracy from $\gm_1$ to
$\gm_m$. Then $(\gm\hskip-.03in.^+,L.^+)$ is an augmented free p-bounded simplicial group
which is
$(G,L_{-1} = L_G)$ in dimension $-1$, and equipped with a word-length metric in non-negative
degrees.

\proclaim{\bf\underbar{Claim 2.1.5}} If the Dehn function of $G$ is polynomial, then there
is a
p-bounded section of sets $s'_1 : \gm^0_0\inj \gm^0_1$ which is a left-inverse to
$\partial_1$
restricted to $\gm^0_1$.
\endproclaim

\prf As $\gm^0_0 = ker(\varepsilon)$ is the subgroup of $\gm_0$ normally generated by the
relators $\Cal W$, an element $\ov w\in \gm^0_0$ may be written
$$
\ov w = \ov w_1^{x_1}\ov w_2^{x_2}\dots \ov w_n^{x_n}
\tag2.1.6
$$
as in (2.1.1), where $x_i\in \gm_0$ and $w_i^{\pm}\in \Cal W\subset \Bbb X_1$.
We use the convention of (2.1.4) to distinguish between $w\in \Cal W$ and
$\ov w = \partial_1(w)\in \gm_0$. A p-bounded section $s'_1$ exists if and only if there are
constants $C,n > 0$ such that for all $w\in \gm_0^0$, there exists a choice of $w_i$ and $x_i$
in (2.1.6) for which
$$
L_1(w_1^{s_0(x_1)}w_2^{s_0(x_2)}\dots w_n^{s_0(x_n)})\le C(1 + L_0(w))^n
$$
Since
$$
L_1(w_1^{s_0(x_1)}w_2^{s_0(x_2)}\dots w_n^{s_0(x_n)})\le \sum_{i=1}^n L_1(w_i) + 2L_0(x_i)
\tag2.1.7
$$
Gersten's Lemma 2.1.4 implies the left-hand side of (2.1.7) is quadratically bounded by
the Dehn function $f_G$ of $G$. Then $f_G$ polynomial implies the claim.\hfill //
\endpf

Continuing with the proof of Theorem 2.1.3, we see that $(\gm\hskip-.03in.^+,L.^+)$ is type
$P(1)$
as defined above. By Theorem A.1 of the appendix, there is an inclusion of augmented
simplicial groups with word-length (and generating sets)
$$
(\gm\hskip-.03in.^+,\Bbb X.^+,L.^+)\hookrightarrow (\wt \gm\hskip-.03in.^+,\wt{\Bbb X}.^+,\wt L.^+)
$$
where $(\wt \gm\hskip-.03in.^+,\wt{\Bbb X}.^+,\wt L.^+)$ is a type $P$ resolution
and $\wt \gm_i = \gm_i$ for $i = 0,1$.
Now set $D^n(\gm\hskip-.03in.,\Bbb Q) = PH^1(\gm_{n-1};\Bbb Q)$ for $n\ge 1$, and $0$ for
$n = 0$. Similarly, let
$E^n(\gm\hskip-.03in.,\Bbb Q) = H^1(\gm_{n-1};\Bbb Q)$ for $n\ge 1$, and $0$ for $n = 0$.
There are coboundary maps
$\delta^n = \sum_{i=1}^n (-1)^i(\partial_i)^*: E^n(\gm\hskip-.03in.,\Bbb Q)
= H^1(\gm_{n-1};\Bbb Q)\to E^{n+1}(\gm\hskip-.03in.,\Bbb Q) = H^1(\gm_{n};\Bbb Q)$
making $(E^*(\gm\hskip-.03in.,\Bbb Q),\delta^*)$ a cocomplex. Because each face map
$\partial_i$ is p-bounded, we also get a well-defined coboundary map
$\delta^n: D^n(\gm\hskip-.03in.,\Bbb Q) =
PH^1(\gm_{n-1};\Bbb Q)\to D^{n+1}(\gm\hskip-.03in.,\Bbb Q)
= PH^1(\gm_{n};\Bbb Q)$ given by the same expression. In addition, since
$(\gm_n,L_n)$ is a free group with word-length metric, there is for each $n\ge 0$
an inclusion $D^n(\gm\hskip-.03in.,\Bbb Q)\hookrightarrow E^n(\gm\hskip-.03in.,\Bbb Q)$
which is clearly compatible with the coboundary maps just defined, yielding an
inclusion of cocomplexes
$$
(D^*(\gm\hskip-.03in.,\Bbb Q),\delta^*)\hookrightarrow(E^*(\gm\hskip-.03in.,\Bbb Q),\delta^*)
\tag2.1.8
$$
By Theorem A.12 of the appendix, there are isomorphisms
$$
\gathered
PH^*(G;\Bbb Q)\cong H^*(D^*(\gm\hskip-.03in.,\Bbb Q),\delta^*)\\
H^*(G;\Bbb Q)\cong H^*(E^*(\gm\hskip-.03in.,\Bbb Q),\delta^*)
\endgathered
\tag2.1.9
$$
under which the inclusion of (2.1.8) induces the transformation
$$
PH^*(G;\Bbb Q)\to H^*(G;\Bbb Q)
\tag2.1.10
$$
Because the generating set for $\gm_i$ is finite for $i = 0,1$, the map in (2.1.8) is
an isomorphism for $* = 1,2$. Together with the injectivity of the map for $*=3$,
Theorem A.12 implies the map in (2.1.10) is an isomorphism for $* = 1,2$.\hfill //
\endpf

If $A$ is a p.s\. $G$-module, it is a p.s\. $\gm_i$-module via the augmentation
$\gm_i\surj G$. One may then replace the coefficient module $\Bbb Q$ by $A$ in the above
discussion. The result is again an inclusion of cocomplexes
$$
(D^*(\gm\hskip-.03in.,A),\delta^*)\hookrightarrow(E^*(\gm\hskip-.03in.,A),\delta^*)
\tag2.1.11
$$
and isomorphisms
$$
\gathered
PH^*(G;A)\cong H^*(D^*(\gm\hskip-.03in.,A),\delta^*)\\
H^*(G;A)\cong H^*(E^*(\gm\hskip-.03in.,A),\delta^*)
\endgathered
\tag2.1.12
$$
under which the inclusion of (2.1.11) induces the transformation
$$
PH^*(G;A)\to H^*(G;A)
\tag2.1.13
$$
The finiteness of $\gm_i$ for $i = 0,1$ implies

\proclaim{\bf\underbar{Corollary 2.1.14}} If $G$ is a finitely-presented group with Dehn
function of polynomial type and $A$ is a p.s\. $G$-module, then
$PH^*(G;A)\cong H^*(G;A)$ for $* = 0,1,2$.
\endproclaim

%%%%%%%%%%%%%%%%%%%%%%%%%%%%%%%%%%%%%%%%%%%%%%%%%

\vskip.3in
\subhead
2.2 Higher Dehn functions and cohomology
\endsubhead
\vskip.2in

Suppose $\pi$ is an $HF^{\infty}$ group, i.e.
one with a classifying space
$B\pi$ the homotopy type of a $CW$ complex with finitely many cells in each dimension.
One can show this is equivalent to the condition that there is a simplicial set $X.$ with
$X_0 = *$ and $|X|\simeq B\pi$ where $X\hskip-.025in.$ has finitely many simplicies in each
dimension. The augmented Kan
loop group $GX\hskip-.025in.^+$ of $X\hskip-.025in.$ is $\pi$ in dimension $-1$, and in
dimension $n\ge 0$ is the free
group on generators $\Bbb X_n = X_{n+1} - s_0(X_n)$, which is a finite set.
We write $GX\hskip-.025in.^+$ as $\gm\hskip-.03in.^+$.
For $n\ge 0$, $L_n$ is the standard word-length function associated with the set of generators
$\Bbb X_n$. The generating set $\Bbb X_0$ determines a generating set for $\pi$, and we take
$L_{-1}$ to be the standard word-length function on $\pi$ associated with those generators.
Then $L.^+ = \{L_n\}_{n\ge -1}$ is an augmented simplicial word-length function, and
$(\gm\hskip-.03in.^+,L.^+)$ is an augmented free simplicial resolution with word-length
function, with $\gm_n$ generated by a finite set $\Bbb X_n$ and $L_n$ the standard
word-length function determined by $\Bbb X_n$. We assume for this section that
$(\gm\hskip-.03in.^+,L.^+)$ is as just
described.
As discussed in the previous section, the fact $\gm\hskip-.03in.^+$ is a resolution
implies the homomorphism
$$
\partial_{n+1}' : \gm^n_{n+1}\to \gm^n_n
\tag2.2.1
$$
induced by the restriction of $\partial_{n+1}$ to $\gm^n_{n+1}$, is surjective for all
$n\ge 0$.
We now define the $n^{th}$ Dehn function $f^{\gm}_n$ associated to $\gm\hskip-.03in.^+$
as the smallest $\Bbb N$-valued function for which there exists a section of sets
$s'_{n+1} : \gm^n_n\to \gm^n_{n+1}$ with
$\partial_{n+1}'\circ s'_{n+1} = id$ and
$$
L_{n+1}(s'_{n+1}(x))\le f^{\gm}_n(L_n(x))\qquad \forall x\in \gm^n_n
\tag2.2.2
$$
The Dehn function of $(\gm\hskip-.03in.^+,L.^+)$ is
$f\hskip-.02in.^{\gm} = \{f^{\gm}_n\}_{n\ge 0}$.
An element $x\in \gm_n^n$ induces a map $\phi_x : S^n\to |\gm.|$ and an element
$y\in\gm_{n+1}^n$ with $\partial_{n+1}(y) = x$ induces a null-homotopy of $\phi_x$. So this
definition is the simplicial analogue of the classical geometric situation where one bounds
the volume of a null homotopy of a map $S^n\to M$ by a function evaluated on the volume of
(the image of) $S^n$. We call two functions $f_1,f_2 : \Bbb N\to \Bbb N$
\underbar{p-equivalent} if there are polynomials $p_1,p_2$ such that $f_1\le p_2\circ f_2$
and $f_2\le p_1\circ f_1$. A long and
tedious argument using simplicial identities proves the following

\proclaim{\bf\underbar{Theorem 2.2.3}} If $f\hskip-.02in.^{\gm}$ is the Dehn function
for $(\gm\hskip-.03in.^+,L.^+)$, then
$\gm(\varepsilon)\hskip-.01in.^+ = ker(\gm.^+\surj \pi)$ admits an extra degeneracy
$s.' = \{s'_{n+1}\}_{n\ge -1}$ satisfying the property that
$$
L_{n+1}(s'_{n+1}(x))\le f'_n(L_n(x))\qquad \forall x\in \gm(\varepsilon)_n
$$
where $f'_n$ is p-equivalent to $f^{\gm}_n$ for all $n$.\hfill\text{//}
\endproclaim

Note that the converse is obvious, since $\gm(\varepsilon)_n^n = \gm_n^n$.
If the Dehn function $f^{\gm}_m$ is polynomial for $m\le n$, then
$(\gm\hskip-.03in.^+,L.^+)$ is a type
$P(n+1)$ resolution, and if it is polynomial for all $m$ then
$(\gm\hskip-.03in.^+,L.^+)$ is type $P$. Finite generation in each degree
implies
$PH^1(\gm_n;A) = H^1(\gm_n;A)$ for $n\ge 0$, where $A$ is an arbitrary p.s\. $\pi$-module.
By the same method as above we have

\proclaim{\bf\underbar{Theorem 2.2.4}} Let $A$ be a p.s\. $\pi$-module.
If $f^{\gm}_m$ is polynomial for
$m\le n$, then
$(\gm\hskip-.03in.^+,L.^+)$ is a type $P(n+1)$ resolution of $\pi$, and the map
$$
PH^*(\pi;A)\to H^*(\pi;A)
$$
is an isomorphism for $*\le n+2$. If $f^{\gm}_m$ is polynomial for all $m$, then $\pi$
satisfies condition (PC1) for all coefficients $A$.\hfill\text{//}
\endproclaim

In the proof of Theorem 2.1.3, we showed how to construct a simplicial group from a
presentation of $\pi$. If that presentation is finite, $\gm_0$ and $\gm_1$ are finitely
generated. If $\pi$ is $HF^{\infty}$, this free simplicial group can be extended to a
resolution $\gm\hskip-.03in.$ of the type used in this section. By Lemma 2.1.4, $f^{\gm}_0$ is p-equivalent to
the Dehn function associated with the presentation. This justifies the term higher Dehn
functions when referring to $\{f^{\gm}_n\}_{n\ge 0}$.

\proclaim{\bf\underbar{Question 2.2.5}} If $f^{\gm}_0$ is polynomial, is $f^{\gm}_n$
polynomial for all $n > 0$?
\endproclaim

A stronger version of the same question is

\proclaim{\bf\underbar{Question 2.2.6}} Is $f^{\gm}_n$ polynomially equivalent to $f^{\gm}_0$
for all $n > 0$?
\endproclaim

We conclude this section with an alternative definition of higher Dehn functions analogous to
that given in [Al1]. We first assume, as before, that $\Bbb Q$ admits a resolution over
$\Bbb Q[\pi]$ which is free and finitely generated in each dimension. This may then be written
as
$$
\Cal R(\pi)_* = \Bbb Q\overset\varepsilon\to\leftarrow \Bbb Q[\pi][S_0]
\overset d_1\to\leftarrow\Bbb Q[\pi][S_1]
\overset d_2\to\leftarrow\dots
\tag2.2.7
$$
where each $S_i$ is a finite set and each differential $d_i$ is a $\Bbb Q[\pi]$-module
homomorphism. Taking the weight function on each set $S_i$ to be identically $1$ gives
each term a p\.f\. $\pi$-module structure. Note also that as each $d_i$ is $\pi$-equivariant
and each $S_i$ is finite, the differentials will be linearly bounded with respect to this
p.f\. $\pi$-module structure. Then $\{f_n\}_{n\ge 0}$ is a sequence of isoperimetric functions
for this resolution if there is a chain contraction
$$
\{s_n:\Cal R (\pi)_{n-1}\rightarrow \Cal R (\pi)_n\}_{n\ge 0}
$$
over $\Bbb Q$ with
$$
|s_{n+1}(a)|\le f_n(|a|)
$$
for all $a\in \Cal R(\pi)_n = \Bbb Q[\pi][S_n]$, for all $n\ge 0$ (where for $x\in \Cal R(\pi)_m$, $|x|$ is the semi-norm of $x$ in the p\.f. module $\Cal R(\pi)_m$ as defined in (1.1.4)). If each $f_n$ is a minimal
isoperimetric function, then it is natural to call the sequence the (higher) Dehn functions
associated to the resolution $\Cal R(\pi)_*$. Because each $S_i$ is finite, there are
equalities
$$
PHom_{\pi}(\Bbb Q[\pi][S_n];A) = Hom_{\pi}(\Bbb Q[\pi][S_n];A)
\tag2.2.8
$$

\proclaim{\bf\underbar{Proposition 2.2.9}} If each of the Dehn functions $\{f_n\}_{n\ge 0}$ is
of polynomial type, there is an isomorphism
$$
PH^*(\pi;A)\overset\cong\to\rightarrow H^*(\pi;A)
$$
\endproclaim

\prf  If the Dehn functions are all of polynomial type, the Comparison Theorem yields an
isomorphism
$$
PH^*(\pi;A)\cong H^*(\{PHom_{\pi}(\Bbb Q[\pi][S_n];A),(d_n)^*\}_{n\ge -1})
$$
The result then follows from the isomorphism in (2.2.8).\hfill //
\endpf

As an application we have

\proclaim{\bf\underbar{Corollary 2.2.10}} If $A$ is a p.s\. $G$-module and $G$ is either finitely generated nilpotent or synchronously combable, then there is an isomorphism
$$
PH^*(\pi;A)\overset\cong\to\rightarrow H^*(\pi;A)
$$
\endproclaim

\prf For nilpotent groups, all the Dehn functions are polynomial, so we may apply the above methods. The result may alternatively be proved by induction on the length of the lower central series and the fact that for abelian central extensions, hypothesis 1.1.1H(k) can be shown to hold for all $k$, allowing for a comparison of Serre spectral sequences. For the second case, we appeal to Gersten's argument in [G5]. Gersten's argument not only shows that $G$ is of type $FP^{\infty}$ (cf. [Al1]), but that the Dehn functions associated to the cellular chain complex of the universal cover, which is of the type in (2.2.7), are all polynomial (not just the first one). The result follows by the above proposition.\hfill //
\endpf

%%%%%%%%%%%%%%%%%%%%%%%%%%%%%%%%%%%%%%%%%%%%%%%%%

\vskip.3in
\subhead
2.3 Linearly and uniformly bounded cohomology
\endsubhead
\vskip.2in

Probably the strongest constraint one can impose on the cochain level
while still retaining enough functoriality
for the homological algebra machinery of section 1.1 is a linear (or Lipschitz) constraint.
Thus, an
\underbar{l-semi-normed $G$-module} (l.s\. $G$-module) is defined as in (1.1.1) except that
in (1.1.1) ii) we replace $n$ by $1$. A homomorphism $f : M\to M'$ of l.s\. $G$-modules is a
$\Bbb Q[G]$-module homomorphism for which there exists $C_1,C_2 > 0$ such that for all
$x'\in \Lambda_{M'}$ there exists $x\in \Lambda_M$ with
$$
\eta_{gx'}(f(a))\le C_1C_2(1+\eta_{gx}(a))
$$
for all $g\in G$ and $a\in M$. As before, the constant $C_1$ may vary with $f$
but is independent of $x'$, while $C_2$ may vary with $x'$ but is independendent of the
other parameters. Note this is slightly more rigid than what one gets when replacing
$n$ by $1$ in (1.1.2). The set of linearly bounded l.s\. $G$-module homomorphisms
from $M$ to $M'$ is denoted by $LHom_G(M,M')$; not requiring $f$ to commute with the
action of $G$ produces the
larger vector space $LHom(M,M')$ on which $G$ acts with fixed-point set $LHom_G(M,M')$.

Finally, an l.f\. resolution of $\Bbb Q$ over $\Bbb Q[G]$ is defined
as in (1.1.9), except that the the differentials are homomorphisms of l.f\. $G$-modules,
and the chain contraction is required to be linearly bounded. Admissible monomorphisms and
epimorphisms are defined in the same manner as before, with linear replacing polynomial. In
this context, Propositions 1.1.5, 1.1.6 and Lemma 1.1.8 carry over to the linear setting.
Moreover, the bar resolution described prior to the Comparison Theorem is an l.f\.
resolution as the reader may easily verify. This provides the resolution for defining
the linearly bounded cohomology of $G$ with coefficients in an l.s\. $G$ module $M$:
$$
LH^*(G;M)\overset def\to = LH^*_G(EG_*;M) = H^*(LHom_G(EG_*,M),\delta^*)
$$
It is not clear at this point if there is a useful Serre spectral sequence in linearly
bounded cohomology (for reasons discussed below). However, the method of proof of the
Comparison Theorem does carry over, yielding

\proclaim{\bf\underbar{Theorem 2.3.1 - Linear Comparison Theorem}} Let $(R_*,d_*)$ be an l.f\.
resolution of $\Bbb Q$ over $\Bbb Q[G]$ and $M$ an l.s\. $G$-module. Then there is an
isomorphism
$$
LH^*_G(EG_*;M)\cong LH^*_G(R_*;M) = H^*(LHom_G(R_*,M),\delta^*)
$$
\endproclaim

Suppose that the resolution $R_*$ satisfies the finiteness condition mentioned in the previous
section; i.e., $R_n = \Bbb Q[G][S_n]$ with $S_n$ finite for each $n\ge 0$. Then there is
an equality
$$
LHom_{G}(\Bbb Q[G][S_n];A) = Hom_{G}(\Bbb Q[G][S_n];A)
\tag2.3.2
$$

It is reasonable to ask under what conditions such a resolution can exist. The answer is:
when $G$ is word-hyperbolic. This is proved by Mineyev in [M3] (if we worked over
$\Bbb Z$ instead of $\Bbb Q$ then the work of Mineyev and others shows that such a resolution
exists \underbar{if and only if} $G$ is word-hyperbolic). Combining Mineyev's results
with the above yields

\proclaim{\bf\underbar{Theorem 2.3.4}} If $G$ is word-hyperbolic, there is an isomorphism
$$
LH^*(G;M)\overset\cong\to\rightarrow H^*(G;M)
$$
for any l.s\. $G$-module $M$.
\endproclaim

In particular, this implies

\proclaim{\bf\underbar{Corollary 2.3.5}} If $G$ is a finitely-generated free group,
equipped with the standard word-length metric, then
$$
LH^*(G;M) = 0
$$
for all $* > 1$.
\endproclaim

This suggests that the linear analogue of Corollary 1.2.17 may hold.

We conclude this section with a short discussion of bounded
cohomology. Because we do not require word-length functions to be proper,
we could define the length function $L_G$ on $G$ by $L(x) = 1$ if $x\ne 1$.
Because this length function is bounded, the $PHom(\,\,\,)$ groups used in 
the computation of $PH^*(G;M)$ are simply those that are uniformly bounded
on basis vectors, yielding an isomorphism
$$
PH^*(G;M)\cong H_b(G;M)
$$
where the right-hand side denotes the bounded cohomology groups of $G$ in the
p.s\. $G$-module $M$.

Any word-length function on $G$ may be realized as the word-length function induced
by a free group equipped with word-length metric $(F,L_F)$ via an appropriate 
surjection $F\surj G$. In fact, there is a universal example of such. Given $(G,L_G)$
let $F$ be the free group on elements $\{1\ne g\in G\}$ and let $L_F$ be the word-length
metric induced by $L_G$, viewed as a weight function on the set $G - \{1\}$. Then
$$
F'\inj F\surj G
$$
is a short-exact sequence of groups with word-length, and so as before
there is an associated Serre spectral and five-term exact sequence. As in Corollary 1.2.18
one has

\proclaim{\bf\underbar{Corollary 2.3.6}} There is a sequence
$$
H^0(G;A_m^*) = H^{0}_b(G;A_m^*)\surj H^{2}_b(G;A_{m-1}^*)\cong H^4_b(G;A_{m-2}^*)\cong\dots\cong  H^{2m}_b(G;\Bbb Q)
\tag2.3.7
$$
where the maps in the sequence occur as differentials in the $E_2^{*,*}$-term of the
appropriate Serre spectral sequence for p-bounded cohomology (and the groups
$A^*_k$ are as defined in section 1.2).
\hfill\text{//}
\endproclaim

In this way one can realize bounded $2m$-dimensional cohomology classes on $G$ as 
$G$-invariant elements of $A^*_m$. This application to bounded cohomology
will be further examined in the sequel to this paper.

%%%%%%%%%%%%%%%%%%%%%%%%%%%%%%%%%%%%%%%%%%%%%%%%%
%%%%%%%%%%%%%%%%%%%%%%%%%%%%%%%%%%%%%%%%%%%%%%%%%

\newpage
\vskip.2in

\centerline{\bf\underbar{Appendix - Type $P$ resolutions}}
\vskip.3in

The following results first appeared in [O1]. We have included them here as they are an essential ingredient in the proofs appearing in section 2.1 and 2.2.
We begin with a demonstration of the existence of type $P$ resolutions.
\vskip.05in

\proclaim{\bf\underbar{Theorem A.1}} Let $(\gm\hskip-.03in .^+, \Bbb X.^+,L.^+)$ be a
triple, where $(\gm\hskip-.03in .^+,L.^+)$ is a p-bounded augmented free simplicial group with
$\pi = \gm_{-1}$. Then there is an inclusion
$$
\iota : \gm\hskip-.03in .^+\hookrightarrow \wt \gm\hskip-.03in .^+
$$
where $(\wt \gm\hskip-.03in .^+,\wt \BbbX.^+,\wt L.^+)$ is a triple and $(\wt \gm\hskip-.03in .^+,\wt L.^+)$
is a type $P$ resolution of $\pi$. Moreover, if $(\gm\hskip-.03in .^+, \Bbb X.^+,L.^+)$ is type $P(m)$, then the construction can be done so that $(\wt \gm_n^+,\wt \BbbX_n^+,\wt L_n^+) = ( \gm_n^+,\BbbX_n^+,L_n^+)$ for $n\le m$.
\endproclaim

\prf We first give the general construction. Denote $(\gm\hskip-.03in .^+,\Bbb X.^+,L.^+)$ by $(\gm(0).^+,\Bbb
X(0).^+,L(0).^+)$.
Note that
$\gm(0)(\varepsilon)_{-1} = \{1\}$, so taking $s'_0 :
\gm(0)(\varepsilon)_{-1}\to \gm(0)(\varepsilon)_0$ as the
inclusion of the trivial group shows that $(\gm(0).^+,\Bbb
X(0).^+,L(0).^+)$
is type $P(0)$.

By induction,  we may assume that a p-bounded free simplicial group
triple\newline
$(\gm(m-1).^+, \BbbX(m-1).^+, L(m-1).^+)$ has been constructed such that
$\gm(m-1)(\varepsilon)\hskip-.015in .^+$ is
$(m-2)$-connected, and admits a contracting degeneracy
$\{s'_{p+1}\}_{0\le p\le m-2}$ through dimension $(m-2)$ which is
p-bounded.

Let $\BbbX(m)'_{m}$ equal the \underbar{set}
$\gm(m-1)(\varepsilon)_{(m-1)} - \{1\}$. For $1\ne
g\in \gm(m-1)(\varepsilon)_{(m-1)}$, we denote by $[g]$ the corresponding
generator in $\BbbX(m)'_{m}$. Let
$$
\gather
\BbbX(m)_j = \BbbX(m-1)_j \,\text{ for }\, j\le (m-1) \;, \\
\BbbX(m)_m = \BbbX(m-1)_m \coprod \BbbX(m)'_{m}\; , \tag{A.2}\\
\BbbX(m)_n = \BbbX(m-1)_n \coprod \{s(\BbbX(m)'_{m})\} \,\text{ for }\, n
> m
\endgather
$$
where the last coproduct is over all iterated degeneracies $s$ from
dimension
$m$ to dimension $n$. Face maps are determined by the following values on
generators
$$
\gather
\partial_j([g]) = s'_{m-1}\partial_j(g)\,
\text{ for }\, 0\le j < m \tag{A.3}\\
\partial_{m}([g]) = g
\endgather
$$

Proceeding as before, we define $\gm(m)_{-1} = \pi$  and $\gm(m)_n$ to be
the free group on $\BbbX(m)_n$ for $n\ge 0$. $L(m)\hskip-.015in .^+$ is uniquely defined
and determined
by the following four properties:
\vskip.2in

(A.4) i) It equals $L(m-1)\hskip-.015in .^+$ on $\BbbX(m-1)\hskip-.015in .^+$
\vskip.05in

(A.4) ii) If $f$ is the proper function on $\BbbX(m)'_{m}$ determined by
$L(m-1)_{m-1}$ restricted to $\gm(m-1)(\varepsilon)_{(m-1)}$, then
$L(m)_m$ is the metric induced by $f$ when restricted to the free group
$F_m$ generated by $\BbbX(m)'_{m}$.
\vskip.05in

(A.4) iii) If $x\in \gm(m)_m$ is written as $x = w_1w_2\dots w_p$ with
$w_{2i-1}$ in $\gm(m-1)_m$ and $w_{2i}\in F_m$, then
$$
L(m)_m(x) = \sum_{i=1}^p L(m)_m(w_i)
$$
\vskip.05in

(A.4) iv) If $s : \gm(m)_m\to \gm(m)_n$ is an iterated degeneracy, then
$L(m)_n(s(x)) = L(m)_m(x)$. If $w = w_1w_2\dots w_q\in \gm(m)_n$ is
a product of degenerate elements $w_{2i}= s(x_i), x_i\in F_m$ and elements $w_{2i-1}$ in
$\gm(m-1)_n$, then $L(m)_n(w) = \sum_{i=1}^q L(m)_n(w_i)$.
\vskip.1in

\noindent That $L(m-1)\hskip-.015in .^+$ is a metric implies
$L(m)\hskip-.015in .^+$ is again a metric
in each non-negative degree, and the contracting degeneracy
$s\hskip-.015in .'$ for $\gm(m)\hskip-.015in .^+$ is now extended through dimension $(m-1)$ as the set
map
$$
\gathered
s'_{m}(1) = 1\\
s'_{m}(g) = [g]\,\text{for } g\ne 1
\endgathered
\tag{A.5}
$$
(A.3) and (A.5) guarantee that $s\hskip-.015in .'$ satisfies the required simplicial
identities
through dimension $(m-1)$. (A.3) -- (A.5) and induction imply that
all of the degeneracy maps (including $s\hskip-.015in .'$) through dimension $(m-1)$ and
all of
the face maps through dimension $m$ are p-bounded. Let
$$
(\wt \gm\hskip-.03in ., \wt \BbbX., \wt L.) =
\varinjlim_m\{(\gm(m).^+,\BbbX(m).^+,L(m).^+)\}\;.
\tag{A.6}
$$
Then $(\wt \gm\hskip-.03in ., \wt {\Bbb X}., \wt L.)$ is a type $P$ resolution. The
inclusion of generating sets $\Bbb X.\hookrightarrow \wt {\Bbb X}.$ induces
the simplicial group monomorphism $\iota : \gm\hskip-.03in .^+\to \wt\gm\hskip-.03in .^+$, which is
the
identity on $\pi = \gm_{-1} = \wt \gm_{-1}$.

Finally if $(\gm\hskip-.03in .^+, \Bbb X.^+,L.^+)$ is type $P(m)$, then in the above sequence we may start with $(\gm(m)\hskip-.01in.^+,\BbbX(m).^+,L(m).^+) = (\gm\hskip-.03in .^+, \Bbb X.^+,L.^+)$ and continue  with the construction by adding generators in simplicial dimensions $n > m$. This verifies the second part of the theorem.\hfill //
\enddemo

\underbar{Example A.7} Let $\pi$ be a countable group equipped with an
$\Bbb N$-valued
word-length function $L$. Let
$\gm\hskip-.03in .^+ = GB\hskip-.025in .\pi^+$ be the augmented
Kan loop group of the non-homogeneous bar construction on $\pi$
(this is the augmented simplicial group associated to the
usual Kan loop group $GB\hskip-.025in .\pi$; cf. [M]). Then the word length function $L$
determines a proper function on the set of $n$-simplicies of
$B\hskip-.025in .\pi$ in
the standard way:
$$
L([g_1,\dots g_n]) = \sum_{i=1}^n L(g_i)
$$
and thus by restriction a proper $\Bbb N^+$-valued function on the
generating set $\Bbb X_{n-1} = B_n\pi - s_0(B_{n-1}\pi)$ of
$(GB\pi)_{n-1}$ for all $n\ge 1$. In non-negative dimensions we then take
$L_n$ to be the metric determined by this proper function. This produces a
resolution $(\gm\hskip-.03in .^+,\Bbb X.^+,L.^+)$ to which we may apply the above
extension theorem. Note also that the word-length functions $L_n$ arising
from this construction are $\Bbb N$-valued, making the word-length
function $\wt L.$ constructed above $\Bbb N$-valued as well.
\vskip.2in

We summarize this as

\proclaim{\bf\underbar{Corollary A.8}} Every countable group $\pi$
admits a type $P$ resolution where the word-length function in
non-negative degrees is an $\Bbb N$-valued metric. Moreover, if $C(\pi)$ is the category
whose objects are p-bounded augmented free simplicial groups equal to
$\pi$ in dimension $-1$ and equipped with word-length metrics in
non-negative degrees, and whose morphisms are p-bounded simplicial group
homomorphisms inducing the identity on $\pi_0$, then the full subcategory
whose objects are type $P$ resolutions is cofinal in $C(\pi)$.
\endproclaim

For the remainder of the section we assume $\gm\hskip-.03in .^+$ is a type $P$ resolution of $G$. We construct a contraction of $\gm\hskip-.03in .^+$ viewed as an 
augmented simplicial set. To begin with, $\gm(\varepsilon)\hskip-.015in .$ admits a 
simplicial contraction
$s\hskip-.015in .' = \{s'_{n+1} : \gm(\varepsilon)_n \to
\gm(\varepsilon)_{n+1}\}_{n\ge 0}$ which is p-bounded.
Now choose a section
$s(0):\gm_{-1}\rightarrowtail \gm_0,\, \varepsilon_0 \circ
s(0) =$ identity, with $s(1) = 1$ and which is minimal with respect to
word-length. Define $s(n) = s^{(n)}_0 \circ s(0): \gm_{-1} \rightarrowtail
\gm_n$.  Note that
$$
\aligned
\varepsilon_n\circ s(n) &= \text{ identity } \qquad \forall\, n\ge 0\; , \\
\partial_i \circ s(n) &= s(n-1) \qquad\quad \forall n\ge 1, 0\le i \le n\;
, \\ s_i \circ s(n-1) &= s(n) \qquad\qquad \forall n\ge 1, 0\le i \le
n-1\; . \endaligned
\tag{A.9}
$$
For $n = -1$, set $\wt s_{n+1} = \wt s_0 = s(0)$.
Note that for arbitrary $g\in \gm_n$, $g(s(n)(\varepsilon_n(g)))^{-1}\in
\gm(\varepsilon)_n$. 
Then when $n\ge 0$
$$
\wt s_{n+1}(g) = s'_{n+1}(g(s(n)(\varepsilon_n(g)))^{-1})s(n+1)(\varepsilon(g))
\tag{A.10}
$$
This defines a map of sets $\wt s_{n+1} :\gm_n\to \gm_{n+1}$. The simplicial identities
imply $\wt s_{*+1} = \{\wt s_{n+1}: \gm_n\to \gm_{n+1}\}_{n\ge -1}$ is a simplicial contraction
of simplicial sets, which by construction is p-bounded for each $n\ge -1$.

Recall from section 2.1 that $D^n(\gm\hskip-.03in .,\Bbb Q) = PH^1(\gm_{n-1};\Bbb Q)$ for $n\ge 1$, 
and $D^n(\gm\hskip-.03in .,\Bbb Q) = 0$ for $n \le 0$. As all face
maps of $\gm\hskip-.03in .$ are p-bounded, there is a well-defined homomorphism
$\delta^n = \sum_{i=0}^{n-1}(-1)^i\partial_i^* : 
D^n(\gm\hskip-.03in .,\Bbb Q)\to D^{n+1}(\gm\hskip-.03in .,\Bbb Q)$, making
$(D^*(\gm\hskip-.03in .,\Bbb Q),\delta^*)$ a cocomplex. Similarly, one defines the cocomplex $(E^*(\gm\hskip-.03in .,\Bbb Q),\delta^*)$ in the same fashion with $H^1$ in place of $PH^1$. As we observed in section 2.1, there is an inclusion of cocomplexes 
$$
(D^*(\gm\hskip-.03in .,\Bbb Q),\delta^*)\hookrightarrow (E^*(\gm\hskip-.03in .,\Bbb Q),\delta^*)
\tag{A.11}
$$

\proclaim{\bf\underbar{Theorem A.12}} For $n\ge 1$ there is an isomorphism of cohomology groups
$$
PH^n(G;\Bbb Q)\cong H^n(D^*(\gm\hskip-.03in .,\Bbb Q),\delta^*)
$$
Moreover, the inclusion of cocomplexes $(D^*(\gm\hskip-.03in .,\Bbb Q),\delta^*)\hookrightarrow (E^*(\gm\hskip-.03in .,\Bbb Q),\delta^*)$ induces, upon passing to cohomology, the transformation $PH^*(G;\Bbb Q)\to H^*(G;\Bbb Q)$.
\endproclaim

\prf Fix an $m\ge 0$ and consider the augmented simplicial abelian group
$$
C(m).^+ = \{[n]\mapsto C_m((B\gm_n);\Bbb Q)\}_{n\ge -1}
$$
The p-bounded contraction $\wt s_{*+1}$ on $\gm\hskip-.03in.^+$ defined above
induces a p-bounded $\Bbb Q$-vector space contraction
on $C(m).^+$ for each $m\ge 0$ given by
$$
B_m\gm_n\ni [g_1,\dots,g_m]\mapsto 
[\wt s_{n+1}(g_1),\dots,\wt s_{n+1}(g_m)]\in B_m\gm_{n+1}\qquad n\ge -1
$$ 
Applying $PHom(\,\,\, ,\Bbb Q)$ to the associated complex
$C(m)_*$, we get a cocontraction above dimension $0$, yielding
$$
\gather
H^n(PHom(C(m)_*,\Bbb Q)) = 0 \,\,\text{for }n > 0\\
H^0(PHom(C(m)_*,\Bbb Q)) = PHom(C_m(BG;\Bbb Q),\Bbb Q)
\endgather
$$
Applying $PHom(\,\,\, ,\Bbb Q)$
to the bi-complex $C_{*,*} = \{C_*(m)\}_{m\ge 0}$ produces a bi-cocomplex. From the computation
of $H^*(PHom(C(m)_*,\Bbb Q))$, we see that filtering by columns produces an $E_1$-term which 
collapses to the cocomplex 
$(PHom(C_*(BG;\Bbb Q),\Bbb Q),\delta^*)$ whose cohomology is 
$PH^*(G;\Bbb Q)$.
Filtering by rows on the other hand yields an $E_1$-term with $E_1^{p,q} = PH^p(\gm_q;\Bbb Q)$.
Now $\gm_q$ is free and $L_q$ is a metric for $q\ge 0$, so by Corollary 1.2.17
$E_1^{p,q} = 0$ for $p > 1$. The
$E_2^{0,*}$-line, which is the cohomology of $(E_1^{0,*},d_1^{0,*})$, is $\Bbb Q$ for $* = 0$ 
and $0$ for $*>0$. There is an isomorphism of cocomplexes
$(D^*(\gm\hskip-.03in .,\Bbb Q),\delta^*)\cong (E_1^{1,*-1},d_1^{1,*-1})$, hence
$$
PH^n(G;\Bbb Q) = E_2^{1,n-1} = 
H^n(D^*(\gm\hskip-.03in .,\Bbb Q),\delta^*)
$$
for all $n\ge 1$. Applying $Hom(\,\,\, ,\Bbb Q)$ in place of $PHom(\,\,\, ,\Bbb Q)$ and repeating the same line of reasoning produces an isomorphism
$$
H^n(G;\Bbb Q) = E_2^{1,n-1} = 
H^n(E^*(\gm\hskip-.03in .,\Bbb Q),\delta^*)
$$
Finally, the natural transformation $PH^*(G;\Bbb Q)\to H^*(G;\Bbb Q)$ in the above context is induced by a map of bicocomplexes coming from the natural transformation $PHom(\,\,\, ,\Bbb Q)\to Hom(\,\,\, ,\Bbb Q)$. On the level of spectral sequences, this induces a map on the $E^{1,*}_1$ line corresponding to the inclusion of (A.11) above, completing the proof. \hfill //
\endpf

%%%%%%%%%%%%%%%%%%%%%%%%%%%%%%%%%%%%%%%%%%%%%%%%%
%%%%%%%%%%%%%%%%%%%%%%%%%%%%%%%%%%%%%%%%%%%%%%%%%

\enddocument